\documentclass[10pt,reqno]{article}
\usepackage{mathrsfs}
\usepackage{times}
\usepackage[T1]{fontenc}
\usepackage{amsmath,amsfonts,amsthm,amssymb,amscd}
\usepackage[dvips]{graphics}
\usepackage{epsfig}

\newcommand{\ncr}[2]{{#1 \choose #2}}
\newtheorem{thm}{Theorem}[section]
\newtheorem{conj}[thm]{Conjecture}

\numberwithin{equation}{section}

\theoremstyle{definition}

\newtheorem{rek}[thm]{Remark}
\newcommand\ben{\begin{enumerate}}
\newcommand\een{\end{enumerate}}
\newcommand\bi{\begin{itemize}}
\newcommand\ei{\end{itemize}}

\newcommand\be{\begin{equation}}
\newcommand\ee{\end{equation}}
\newcommand\bea{\begin{eqnarray}}
\newcommand\eea{\end{eqnarray}}

\newcommand{\gep}{\epsilon}  %lowercase epsilon
\newcommand{\gl}{\lambda}    %lowercase lambda
\newcommand{\wglt}{\widetilde{\lambda_\pm}}

\numberwithin{equation}{section}

\begin{document}

\title{The Distribution of the Largest Non-trivial
Eigenvalues in Families of Random Regular Graphs}

%\subjclass[2000]{ 05C80 (primary), 05C50, 15A52 (secondary). }
%\keywords{Ramanujan Graphs, Random Graphs, Second Largest
%Eigenvalue, Tracy-Widom Distribution}

\author{Steven J. Miller\thanks{Department of Mathematics, Brown
University, 151 Thayer Street, Providence, RI 02912, E-mail:
\texttt{Steven.Miller.MC.96@aya.yale.edu, sjmiller@math.brown.edu}. \emph{Future address:} Department of Mathematics, Williams College, Williamstown, MA 01267},\ \ \
Tim Novikoff\thanks{Center for Applied Math, Cornell University, Ithaca, NY 14853, E-mail: \texttt{tnovikoff@gmail.com}}\ \ \ and Anthony Sabelli
\thanks{Department of Mathematics, Brown University, Providence, RI 02912, E-mail: \texttt{Anthony\underline{\ }Sabelli@brown.edu}}}

\date{Received November 22\textsuperscript{th}, 2006,
revised February \textsuperscript{th}, 2008.}

%\begin{document}

\maketitle
%\centerline{Mathematics Department}
% \centerline{Brown University}
% \centerline{151 Thayer Street}
% \centerline{Providence, RI 02912}

\begin{center} 2000 AMS Subject
Classification: 05C80 (primary), 05C50, 15A52 (secondary).\end{center}

\begin{abstract} Recently Friedman proved Alon's conjecture for many
families of $d$-regular graphs, namely that given any $\epsilon
> 0$ ``most'' graphs have their largest non-trivial eigenvalue at most
$2\sqrt{d-1}+\gep$ in absolute value; if the absolute value of the largest non-trivial eigenvalue is at most
$2\sqrt{d-1}$ then the graph is said to be Ramanujan. These graphs
have important applications in communication network theory,
allowing the construction of superconcentrators and nonblocking
networks, coding theory and cryptography. As many of these
applications depend on the size of the largest non-trivial positive and negative eigenvalues, it
is natural to investigate their distributions. We show these are
well-modeled by the $\beta=1$ Tracy-Widom distribution for several families. If the
observed growth rates of the mean and standard deviation as a
function of the number of vertices holds in the limit, then in the
limit approximately 52\% of $d$-regular graphs from bipartite families should be Ramanujan, and about $27\%$ from non-bipartite families (assuming the largest positive and negative eigenvalues are independent).
\end{abstract}

%\tableofcontents

\section{Introduction}

\subsection{Families of Graphs}

In this paper we investigate the distribution of the largest
non-trivial eigenvalues associated to $d$-regular undirected graphs\footnote{An
undirected graph $G$ is a collection of vertices $V$ and edges $E$
connecting pairs of vertices. $G$ is simple if there are no multiple
edges between vertices, $G$ has a self-loop if a vertex is connected
to itself, and $G$ is connected if given any two vertices $u$ and
$w$ there is a sequence of vertices $v_1, \dots, v_n$ such that
there is an edge from $v_i$ to $v_{i+1}$ for $i \in \{0, \dots,
n+1\}$ (where $v_0 = u$ and $v_{n+1} = w$).}. A graph $G$ is
bipartite if the vertex set of $G$ can be split into two disjoint
sets $A$ and $B$ such that every edge connects a vertex in $A$ with
one in $B$, and $G$ is $d$-regular if every vertex is connected to
exactly $d$ vertices. To any graph $G$ we may associate a real
symmetric matrix, called its adjacency matrix, by setting $a_{ij}$
to be the number of edges connecting vertices $i$ and $j$. Let us
write the eigenvalues of $G$ by $\gl_1(G) \ge \cdots \ge \gl_N(G)$,
where $G$ has $N$ vertices. We call any eigenvalue equal to $\pm d$ a trivial eigenvalue (there is an eigenvalue of $-d$ if and only if the graph is bipartite), and all other eigenvalues are called non-trivial.

The eigenvalues of the adjacency matrix provide much information about the graph. We give two such
properties to motivate investigations of the eigenvalues; see
\cite{DSV,Sar1} for more details.

First, if $G$ is $d$-regular then $\gl_1(G) = d$ (the corresponding
eigenvector is all $1$'s); further, $\gl_2(G) < d$ if and only if
$G$ is connected. Thus if we think of our graph as a network,
$\gl_2(G)$ tells us whether or not all nodes can communicate with
each other. For network purposes, it is natural to restrict to
connected graphs without self-loops.

Second, a fundamental problem is to construct a well-connected
network so that each node can communicate with any other node
``quickly'' (i.e., there is a short path of edges connecting any two
vertices). While a simple solution is to take the complete graph as
our network, these graphs are expensive: there are $N$ vertices and
$\ncr{N}{2} = N(N-1)/2$ edges. We want a well-connected network where the number of edges grows linearly with
$N$. Let $V$ be the set of vertices for a graph $G$, and $E$ its set
of edges. The boundary  $\partial{U}$ of a $U \subseteq V$ is the
set of edges connecting $U$ to $V\setminus U$. The expanding
constant $h(G)$ is
\begin{equation}
h(G) \ :=\ \inf \left\{{\frac{|\partial{U}|}{\min(|U|,|V\setminus
U|)} : U \subset V,\ |U| > 0}\right\},
\end{equation} and measures the connectivity of $G$. If
$\{G_m\}$ is a family of connected $d$-regular graphs, then we call
$\{G_m\}$ a family of expanders if $\lim_{m\to\infty}|G_m| =\infty$
and there exists an $\epsilon > 0$ such that for all $m$, $h(G_m)
\ge \epsilon$. Expanders have two very important properties: they
are sparse ($|E|$ grows at most linearly with $|V|$), and they are
highly connected (the expanding constants have a positive lower
bound). These graphs have important applications in communication
network theory, allowing the construction of superconcentrators and
nonblocking networks \cite{Bien,Pi}, as well as applications to
coding theory \cite{SS} and cryptography \cite{GILVZ}; see \cite{Sar2} for a brief introduction to expanders. The
Cheeger-Buser inequalities\footnote{The name is from an analogy with
the isoperimetric constant of a compact Riemann manifold.} (due to
Alon and Milman \cite{AM} and Dodziuk \cite{Do}) give upper and
lower bounds for the expanding constant of a finite $d$-regular
connected graph in terms of the spectral gap (the separation between
the first and second largest eigenvalues) $d-\lambda_2(G)$:
\begin{equation}
\frac{d-\lambda_2(G)}{2}\ \le\ h(G)\ \le\
2\sqrt{2d(d-\lambda_2(G))}.
\end{equation}
Thus we have a family of expanders if and only if there exists an
$\epsilon > 0$ such that for all $m$, $d-\lambda_2(G_m) \ge
\epsilon$. Finding graphs with small $\lambda_2(G)$ lead to large
spectral gaps and thus sparse, highly connected graphs.
%; however, $d-\lambda_2(G_m)$ cannot be too large.

For many problems, the behavior is controlled by the largest absolute value of a non-trivial eigenvalue. We write $\lambda_+(G)$ (resp., $\lambda_-(G)$) for the largest non-trivial positive eigenvalue (resp., the most negative non-trivial eigenvalue) of $G$, and set $\lambda(G) = \max\left(|\lambda_+(G)|,|\lambda_-(G)|\right)$. Alon-Boppana, Burger, and
Serre proved that for any family $\{G_m\}$ of finite connected
$d$-regular graphs with $\lim_{m\to\infty} |G_m| = \infty$, we have
$\liminf_{m \rightarrow \infty} \lambda(G_m)\ \ge\ 2 \sqrt{d-1}$; in fact, Friedman \cite{Fr1} proved that if $G$ is a $d$-regular ($d\ge 3$) graph with $n$ vertices, then \be\label{eq:friedmanlowerbound} \lambda(G) \ \ge \ 2\sqrt{d-1} \cdot \left(1- \frac{2\pi^2}{(\log_{d-1} n)^2} + O\left(\frac{1}{(\log_{d-1} n)^4}\right)\right). \ee
Thus we are led to search for graphs with $\gl(G) \le
2\sqrt{d-1}$; such graphs are called Ramanujan\footnote{Lubotzky,
Phillips and Sarnak \cite{LPS} construct an infinite family of $(p +
1)$-regular Ramanujan graphs for primes $p \equiv 1\bmod 4$. Their
proof uses the Ramanujan conjecture for bounds on Fourier
coefficients of cusp forms, which led to the name Ramanujan
graphs.} (see \cite{Mur} for a nice survey). Explicit constructions are known when $d$ is $3$
\cite{Chiu} or $q+1$, where $q$ is either an odd prime
\cite{LPS,Mar} or a prime power \cite{Mor}.

Alon \cite{Al} conjectured that as $N\to\infty$, for $d\ge 3$ and
any $\gep > 0$, ``most'' $d$-regular graphs on $N$ vertices have
$\gl(G) \le 2 \sqrt{d-1} + \gep$; it is known that the
$2\sqrt{d-1}$ cannot be improved upon. Upper bounds on
$\lambda(G)$ of this form give a good spectral gap. Recently,
Friedman \cite{Fr2} proved Alon's conjecture for many models of
$d$-regular graphs. Our goal in this work is to numerically
investigate the distribution of $\gl_\pm(G)$ and $\gl(G)$ for these and other
families of $d$-regular graphs. By identifying the limiting
distribution of these eigenvalues, we are led to the
conjecture that for many families of $d$-regular graphs, in the
limit as the number of vertices tends to infinity the probability a
graph in the family has $\lambda(G) \le 2\sqrt{d-1}$ tends to
approximately 52\% if the family is bipartite, and about $27\%$ otherwise.

Specifically, consider a family $\mathcal{F}_{N,d}$ of $d$-regular
graphs on $N$ vertices. For each $G\in\mathcal{F}_{N,d}$, we study
\be\label{eq:norm2ndlargestevfam} \widetilde{\gl_\pm}(G) \ = \
\frac{\left|\gl_\pm(G)\right| - 2\sqrt{d-1}+c_{\mu,N,d,\pm}
N^{m_\pm(\mathcal{F}_{N,d})}}{c_{\sigma,N,d,\pm} N^{s_\pm(\mathcal{F}_{N,d})}};
\ee we use $m$ for the first exponent as it arises from studying the
means, and $s$ for the second as it arises from studying the
standard deviations. Our objective is to see if, as $G$ varies in a
family $\mathcal{F}_{N,d}$, whether or not $\wglt(G)$ converges to a
universal distribution as $N\to\infty$. We therefore subtract off
the sample mean and divide by the standard deviation to obtain a
mean 0, variance 1 data set, which will facilitate comparisons to
candidate distributions. We write the subtracted mean as a sum of
two terms. The first is $2\sqrt{d-1}$, the expected mean as
$N\to\infty$. The second is the remaining effect, which is observed to
be negative (see the concluding remarks in \cite{Fr2} and \cite{HLW}), and is found
to be negative in all our experiments. We shall assume in our
discussions below that $c_{\mu,N,d,\pm}<0$. Of particular interest is
whether or not $m_\pm(\mathcal{F}_{N,d})-s_\pm(\mathcal{F}_{N,d}) < 0$. If this is negative (for both $\lambda_\pm(G)$), if $\wglt(G)$ converges to a universal
distribution, and if $\lambda_+(G)$ and $\lambda_-(G)$ are independent for the non-bipartite families, then in the limit a positive percent of graphs in
$\mathcal{F}_{N,d}$ are \emph{not} Ramanujan. This follows from the
fact that, for $|\lambda_\pm(G)|$, in the limit a negligible fraction of the standard
deviation suffices to move beyond $2\sqrt{d-1}$; if
$m_\pm(\mathcal{F}_{N,d})-s_\pm(\mathcal{F}_{N,d})
> 0$ then we may move many multiples of the standard deviation and
still be below $2\sqrt{d-1}$ (see Remark \ref{rek:whysmdiffmatters}
for a more detailed explanation).

\begin{rek}[Families of $d$-regular
graphs]\label{rek:familiesdreggraphs} We describe the families we
investigate. For convenience in our studies we always take $N$ to be
even. Friedman \cite{Fr2} showed that for fixed $\gep$, for the
families $\mathcal{G}_{N,d}$, $\mathcal{H}_{N,d}$ and
$\mathcal{I}_{N,d}$ defined below, as $N\to\infty$ ``most''
graphs\footnote{Friedman shows that, given an $\epsilon > 0$, with
probability at least $1 - c_{\mathcal{F}_d}N^{-\tau(\mathcal{F}_d)}$
we have $\lambda(G) \le 2\sqrt{d-1}+\gep$ for
$G\in\mathcal{F}_{N,d}$, and with probability at least
$\widetilde{c}_{\mathcal{F}_d}N^{-\widetilde{\tau}(\mathcal{F}_d)}$
we have $\lambda(G) > 2\sqrt{d-1}$; see \cite{Fr2} for the values
of the exponents.} have $\gl(G) \le 2 \sqrt{d-1} + \gep$.

\bi

\item \emph{$\mathcal{B}_{N,d}$.} We let $\mathcal{B}_{N,d}$
denote the set of $d$-regular bipartite graphs on $N$ vertices. We
may model these by letting $\pi_1$ denote the identity permutation
and choosing $d-1$ independent permutations of $\{1,\dots,N/2\}$.
For each choice we consider the graph with edge set \be E: \ \left\{
(i,\pi_j(i)+ N/2): i \in \{1,\dots, N/2\}, j \in \{1,\dots,
d\}\right\}.\ee

\item \emph{$\mathcal{G}_{N,d}$.} For $d$ even, let $\pi_1, \dots, \pi_{d/2}$
be chosen independently from the $N!$ permutations of
$\{1,\dots,N\}$.  For each choice of $\pi_1,\dots,\pi_{d/2}$ form
the graph with edge set \be E: \ \left\{(i,\pi_j(i)), \
(i,\pi_j^{-1}(i)): i \in \{1,\dots,N\}, j \in
\{1,\dots,d/2\}\right\}. \ee Note $\mathcal{G}_{N,d}$ can have
multiple edges and self-loops, and a self-loop at vertex $i$
contribute $2$ to $a_{ii}$.

\item \emph{$\mathcal{H}_{N,d}$.} These are constructed in the same
manner as $\mathcal{G}_{N,d}$, with the additional constraint that
the permutations are chosen independently from the $(N-1)!$
permutations whose cyclic decomposition is one cycle of length $N$.

\item \emph{$\mathcal{I}_{N,d}$.} These are constructed similarly, except
instead of choosing $d/2$ permutations we choose $d$ perfect
matchings; the $d$ matchings are independently chosen from the
$(N-1)!!$ perfect matchings.\footnote{For example, if $d=3$ and
$N=8$, our three permutations might be $(43876152)$, $(31248675)$
and $(87641325)$. Each permutation generates $8/2 = 4$ edges. Thus
the first permutation gives edges between vertices $4$ and $3$,
between $8$ and $7$, between $6$ and $1$, and between $5$ and $2$. A
permutation whose cyclic decomposition is one cycle of length $N$
can be written $N$ different ways (depending on which element is
listed first). This permutation generates two different perfect
matchings, depending on where we start. Note there are no
self-loops.}

\item \emph{Connected and Simple Graphs.} If $\mathcal{F}_{N,d}$ is
any of the families above ($\mathcal{B}_{N,d}$, $\mathcal{G}_{N,d}$,
$\mathcal{H}_{N,d}$ or $\mathcal{I}_{N,d}$), let
$\mathcal{CF}_{N,d}$ denote the subset of graphs that are connected
and $\mathcal{SCF}_{N,d}$ the subset of graphs that are simple and
connected.

\ei
\end{rek}

\begin{rek}\label{rek:bipartiteevaluessymm} The eigenvalues of bipartite graphs are symmetric about zero. We sketch the proof. Let $G$ be a bipartite graph with $2N$ vertices. Its adjacency matrix is of the form $A(G)=\left({Z \atop B} {B \atop Z}\right)$, where $Z$ is the $N\times N$ zero matrix and $B$ is an $N \times N$ matrix. Let $J =\left( {Z \atop -I} {I \atop Z}\right)$ where $I$ is the $N\times N$ identity matrix. Simple calculations show $J^{-1} = -J$ and $J^{-1} A(G) J = -A(G)$. Noting similar matrices have the same eigenvalues, we see the eigenvalues of $A(G)$ must be symmetric about zero. \end{rek}

\subsection{Tracy-Widom Distributions}\label{subsec:TWdistr}

We investigate in detail the distribution of $\lambda_\pm(G)$ for
$d$-regular graphs related to two of the families above, the perfect
matching family $\mathcal{I}_{N,d}$ and the bipartite family
$\mathcal{B}_{N,d}$ (by Remark \ref{rek:bipartiteevaluessymm} we need only study $\lambda_+(G)$ for the bipartite family). Explicitly, for $N$ even we study
$\mathcal{CI}_{N,d}$, $\mathcal{SCI}_{N,d}$, $\mathcal{CB}_{N,d}$,
and $\mathcal{SCB}_{N,d}$; we restrict to connected graphs as $d$ is a multiple eigenvalue for disconnected graphs. As
$d$ and $N$ increase, so too does the time required to uniformly
choose a simple connected graph from our families; we concentrate on
$d \in \{3,4,7,10\}$ and $N \le 20000$. As there are known
constructions of Ramanujan graphs for $d$ equal to $3$ or $q+1$
(where $q$ is either an odd prime or a prime power), $d=7$ is the
first instance where there is no known explicit construction to
produce Ramanujan graphs. In the interest of space we report in detail on the $d=3$ computations for $\gl_+(G)$. We remark briefly on the other computations and results, which are similar and are available upon request from the authors; much of the data and programs used are available at \begin{center} \texttt{http://www.math.princeton.edu/mathlab/ramanujan/} \end{center}

\begin{figure}
\begin{center}
\caption{\label{fig:plotsTW124} Plots of the three Tracy-Widom
distributions: $f_1(s)$ has the smallest maximum amplitude, then $f_2(s)$ and then $f_4(s)$.} \scalebox{1}{\includegraphics{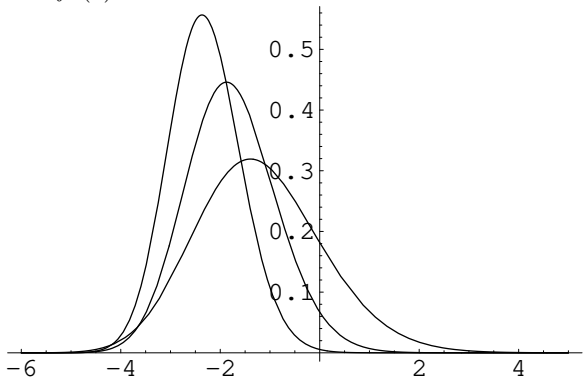}}
\end{center}\end{figure}

We conjecture that the distributions of $\gl_\pm(G)$ are independent in non-bipartite families and each converges to the $\beta=1$ Tracy-Widom distribution (see Conjecture \ref{conj:32} for exact statements). We
summarize our numerical investigations supporting this conjecture in
\S\ref{subsec:summexpresconj}, and content ourselves here with
describing why it is natural to expect the $\beta=1$ Tracy-Widom
distribution to be the answer. The Tracy-Widom distributions model
the limiting distribution of the normalized largest eigenvalues for
many ensembles of matrices. There are three distributions
$f_\beta(s)$: (i) $\beta = 1$, corresponding to orthogonal symmetry
(GOE); (ii) $\beta = 2$, corresponding to unitary symmetry (GUE);
(iii) $\beta = 4$, corresponding to symplectic symmetry (GSE). These
distributions can be expressed in terms of a particular Painlev\'e
II function, and are plotted in Figure \ref{fig:plotsTW124}.

We describe some of the problems where the Tracy-Widom distributions
arise, and why the $\beta = 1$ distribution should describe the
distributions of $\gl_\pm(G)$. The first is in the
distribution of the largest eigenvalue (as $N\to\infty$) in the
$N\times N$ Gaussian Orthogonal, Unitary and Symplectic Ensembles
\cite{TW2}. For example, consider the $N\times N$ Gaussian
Orthogonal Ensemble. From the scaling in Wigner's Semi-Circle Law
\cite{Meh2,Wig}, we expect the eigenvalues to be of order
$\sqrt{N}$. Denoting the largest eigenvalue of $A$ by $\lambda_{{\rm
max}}(A)$, the normalized largest eigenvalue
$\widetilde{\lambda}_{{\rm max}}(A)$ satisfies
\be\label{eq:normlargevTW1gauss} \lambda_{{\rm max}}(A) \ = \
2\sigma \sqrt{N} +\frac{\widetilde{\lambda}_{{\rm
max}}(A)}{N^{1/6}}; \ee here $\sigma$ is the standard
deviation of the Gaussian distribution of the off-diagonal entries,
and is often taken to be $1$ or $1/\sqrt{2}$. As $N\to\infty$ the
distribution of $\widetilde{\lambda}_{{\rm max}}(A)$ converges to
$f_1(s)$. The Tracy-Widom distributions also arise in combinatorics
in the analysis of the length of the largest increasing subsequence
of a random permutation and the number of boxes in rows of random
standard Young tableaux \cite{BDJ,BOO,BR1,BR2,Jo1}, in growth
problems \cite{BR3,GTW,Jo3,PS1,PS2}, random tilings \cite{Jo2}, the
largest principal component of covariances matrices \cite{So},
queuing theory \cite{Ba,GTW}, and superconductors \cite{VBAB}; see
\cite{TW3} for more details and references.

It is reasonable to conjecture that, appropriately normalized, the
limiting distributions of $\gl_\pm(G)$ in the
families of $d$-regular graphs considered by Friedman converges to
the $\beta=1$ Tracy-Widom distribution (the largest eigenvalue is
always $d$). One reason for this is that to any graph $G$ we may
associate its adjacency matrix $A(G)$, where $a_{ij}$ is the number
of edges connecting vertices $i$ and $j$. Thus a family of
$d$-regular graphs on $N$ vertices gives us a sub-family of $N\times
N$ real symmetric matrices, and real symmetric matrices typically
have $\beta = 1$ symmetries. While McKay \cite{McK} showed that for
fixed $d$ the density of normalized eigenvalues is different than
the semi-circle found for the GOE (though as $d\to\infty$ the
limiting distribution does converge to the semi-circle), Jakobson,
Miller, Rivin and Rudnick \cite{JMRR} experimentally found that the
spacings between adjacent normalized eigenvalues agreed with the
GOE. As the spacings in the bulk agree in the limit, it is plausible
to conjecture that the spacings at the edge agree in the limit as
well; in particular, that the density of the normalized second
largest eigenvalue converges to $f_1(s)$.

\subsection{Summary of Experiments, Results and
Conjectures}\label{subsec:summexpresconj}

We numerically investigated the eigenvalues for the
families $\mathcal{CI}_{N,d}$, $\mathcal{SCI}_{N,d}$,
$\mathcal{CB}_{N,d}$ and  $\mathcal{SCB}_{N,d}$. Most of the
simulations were performed on a 1.6GHz Centrino processor running
version 7 of Matlab over several months; the data indicates that the
rate of convergence is probably controlled by the logarithm of the
number of vertices, and thus there would not be significant gains in
seeing the limiting behavior by switching to more powerful systems.\footnote{In fact, many quantities and results related to these families of graphs are controlled by the logarithm of the number of vertices. For example, a family of graphs is said to have large girth if the girths are greater than a constant times the logarithm of the number of vertices (see page 10 of \cite{DSV}). For another example, see \eqref{eq:friedmanlowerbound}.} The data is available at \begin{center} \texttt{http://www.math.princeton.edu/mathlab/ramanujan/} \end{center}

We varied $N$ from $26$ up to $50,000$. For each $N$ we randomly
chose $1000$ graphs $G$ from the various ensembles, and calculated $\lambda_\pm(G)$. Letting
$\mu_{\mathcal{F}_{N,d,\pm}}^{{\rm sample}}$ and
$\sigma_{\mathcal{F}_{N,d,\pm}}^{{\rm sample}}$ denote the mean and
standard deviation of the sample data (these are functions of $N$,
and $\lim_{N\to\infty} \mu_{\mathcal{F}_{N,d,\pm}}^{{\rm sample}} =
2\sqrt{d-1}$), we studied the distribution of
\be\label{eq:norm2ndevsampledata} \left(\lambda_\pm(G) -
\mu_{\mathcal{F}_{N,d,\pm}}^{{\rm sample}}\right)\ \Big/\
\sigma_{\mathcal{F}_{N,d,\pm}}^{{\rm sample}}. \ee This normalizes our
data to have mean 0 and variance 1, which we compared to the $\beta
= 1$ Tracy-Widom distribution; as an additional test, we also compared our data to $\beta = 2$ and $4$ Tracy-Widom distributions, as well as the standard normal.

\begin{figure}
\begin{center}
\caption{\label{fig:plotsTWnorm124} Plots of the three Tracy-Widom
distributions, normalized to have mean 0 and variance 1, and the
standard normal.} \scalebox{1}{\includegraphics{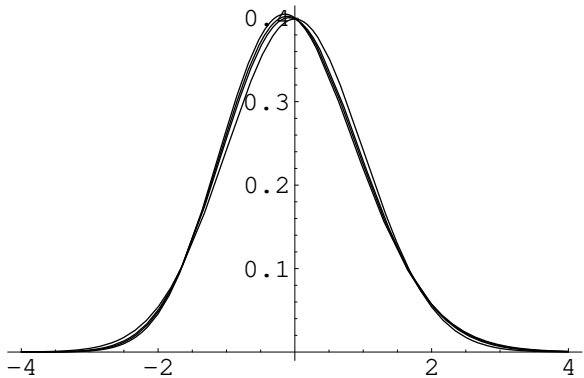}}
\end{center}\end{figure}

Before stating our results, we comment on some of the difficulties
of these numerical investigations.\footnote{Another difficulty is that the Matlab code was originally written to investigate bipartite graphs. The symmetry of the eigenvalues allowed us to just look at the second largest eigenvalue; when we ran the code for non-bipartite graphs, we originally did not realize this had been hardwired. Thus there we were implicitly assuming $\lambda(G)=\lambda_+(G)$, which is frequently false for non-bipartite graphs. This error led us to initially conjecture 52\% of these graphs are Ramanujan in the limit, instead of the 27\% we discuss later.} If $g(s)$ is a probability
distribution with mean $\mu$ and variance $\sigma^2$, then $\sigma
g(\sigma x + \mu)$ has mean 0 and variance 1. As we do not know the
normalization constants in \eqref{eq:norm2ndlargestevfam} for the
second largest eigenvalue, it is natural to study
\eqref{eq:norm2ndevsampledata} and compare our sample distributions
to the normalized $\beta = 1$ Tracy-Widom distribution\footnote{The
Tracy-Widom distributions \cite{TW1} could have been defined in an
alternate way as mean zero distributions if lower order terms had
been subtracted off; as these terms were kept, the resulting
distributions have non-zero means. These correction factors vanish
in the limit, but for finite $N$ result in an $N$-dependent
correction (we divide by a quantity with the same $N$-dependence, so
the resulting answer is a non-zero mean). This is similar to other
situations in number theory and random matrix theory. For example,
originally ``high'' critical zeros of $\zeta(s)$ were shown to be
well-modeled by the $N\to\infty$ scaling limits the $N\times N$ GUE
ensemble \cite{Od1,Od2}; however, for zeros with imaginary part
about $T$ a better fit is obtained by using finite $N$ (in
particular, $N \sim \log T$; see \cite{KeSn}).}. In fact, even if we
did know the constants it is still worth normalizing our data in
order to determine if other distributions, appropriately scaled,
provide good fits as well. As remarked in \S\ref{subsec:TWdistr},
there are natural reasons to suspect that the $\beta = 1$
Tracy-Widom is the limiting distribution; however, as Figure
\ref{fig:plotsTWnorm124} shows, if we normalize the three
Tracy-Widom distributions to have mean 0 and variance 1 then they
are all extremely close to the standard normal. The fact that
several different distributions can provide good fits to the data is
common in random matrix theory. For example, Wigner's
surmise\footnote{Wigner conjectured that as $N\to\infty$ the spacing
between adjacent normalized eigenvalues in the bulk of the spectrum
of the $N\times N$ GOE ensemble tends to $p_W(s) = (\pi
s/2)\exp\left(-\pi s^2/4\right)$. He was led to this by assuming:
(1) given an eigenvalue at $x$, the probability that another one
lies $s$ units to its right is proportional to $s$; (2) given an
eigenvalue at $x$ and $I_1,I_2,I_3,\dots$ any disjoint intervals to
the right of $x$, then the events of observing an eigenvalue in
$I_j$ are independent for all $j$; (3) the mean spacing between
consecutive eigenvalues is $1$.} for the spacings between adjacent
normalized eigenvalues in the bulk of the spectrum is extremely
close to the actual answer (and in fact Wigner's surmise is often
used for comparison purposes, as it is easier to plot than the
actual answer\footnote{The distribution is $(\pi^2/4) d^2
\Psi/dt^2$, where $\Psi(t)$ is (up to constants) the Fredholm
determinant of the operator $f\to\int_{-t}^tK*f$ with kernel
$K=\frac{1}{2\pi}\left(\frac{\sin (\xi-\eta)}{\xi-\eta}+\frac{\sin
(\xi+\eta)}{\xi+\eta}\right)$.}). While the two distributions are
quite close (see \cite{Gau,Meh1,Meh2}) and both often provide good
fits to data, they are unequal and it is the Fredholm determinant
that is correct\footnote{While this is true for number-theoretic
systems with large numbers of data points, there is often not enough data
for physical systems to make a similar claim. The number of energy levels from heavy nuclei in nuclear
physics is typically between 100 and 2000, which can be insufficient to
distinguish between GOE and GUE behavior (while we expect GOE from
physical symmetries, there is a maximum of about a 2\% difference in
their cumulative distribution functions). Current research in
quantum dots (see \cite{Alh}) shows promise for obtaining
sufficiently large data sets to detect such subtle differences.}. We
see a similar phenomenon, as for many of our data sets we obtain
good fits from the three normalized Tracy-Widom distributions and
the standard normal. It is therefore essential that we find a
statistic sensitive to the subtle differences between the four
normalized distributions.

We record the mean, standard deviation, and the percent of the mass
to the left of the mean for the three Tracy-Widom distributions (and
the standard normal) in Table \ref{table:paramsTW}. The fact that
the four distributions have different percentages of their mass to
the left of the mean gives us a statistical test to determine which
of the four distributions best models the observed data.

\begin{table}[ht]
\begin{center}
\caption{\label{table:paramsTW} Parameters for the Tracy-Widom
distributions (before being normalized to have mean 0 and variance 1). $F_\beta$ is the cumulative distribution function for
$f_\beta$, and $F_\beta(\mu_\beta)$ is the mass of $f_\beta$ to
the left of its mean.}
\begin{tabular}{|cccc|}
  \hline
 \ \ \ \  \ \ \ \ &\ \ \ \ Mean $\mu$ \ \ \ \ & Standard Deviation $\sigma$ &
  \ \ \ \ $F_\beta(\mu_\beta)$\ \ \ \  \\
  \hline
TW($\beta=1$) & -1.2065 & 1.26798 & 0.519652\\
TW($\beta=2$) & -1.7711 & 0.90177 & 0.515016\\
TW($\beta=4$) & -2.3069 & 0.71953 & 0.511072\\
Standard Normal &\ 0.0000 & 1.00000 & 0.500000\\ \hline
\end{tabular}
\end{center}
\end{table}

Thus, in addition to comparing the distribution of the normalized
eigenvalues in \eqref{eq:norm2ndevsampledata} to the normalized
Tracy-Widom distributions, we also computed the percentage of time
that $\gl_\pm(G)$ was less than the sample mean. We
compared this percentage to the three different values for the
Tracy-Widom distribution and the value for the standard normal
(which is just .5). As the four percentages are different, this
comparison provides evidence that, of the four distributions, the
second largest eigenvalues are modeled \emph{only} by a $\beta=1$
Tracy-Widom distribution.

%\bigskip \bigskip
%\textbf{Actually, it might be interesting to investigate the
%distribution of $\gl_2(G)$ when $G$ runs through one of the known
%Ramanujan families as well.}\\

We now briefly summarize our results and the conjecture they
suggest. We concentrate on the families (see Remark
\ref{rek:familiesdreggraphs} for definitions) $\mathcal{CI}_{N,d}$,
$\mathcal{SCI}_{N,d}$, $\mathcal{CB}_{N,d}$ and
$\mathcal{SCB}_{N,d}$ with $d\in \{3,4\}$, as well as
$\mathcal{CI}_{N,7}$ and $\mathcal{CI}_{N,10}$. For each $N\in\{$26,
32, 40, 50, 64, 80, 100, 126, 158, 200, 252, 316, 400, 502, 632,
796, 1002, 1262, 1588, 2000, 2516, 3168, 3990, 5022, 6324, 7962,
10022, 12618, 15886, 20000$\}$, we randomly chose 1000 graphs from
each family. We analyze the data for the $3$-regular
graphs in \S\ref{sec:3regstuff}. As the results are similar, the data and analysis for the other families are available online at \texttt{http://www.math.princeton.edu/mathlab/ramanujan/} (where we include our data for $d=3$ as well).\\

\bi

\item \textbf{$\chi^2$-tests for goodness of fit.}
$\chi^2$-tests show that the distribution of the normalized
eigenvalues $\gl_\pm(G)$ are well modeled by a $\beta = 1$ Tracy-Widom
distribution, although the other two Tracy-Widom distributions and
the standard normal also provide good fits; see Tables
\ref{table:chisquare3reg} and \ref{table:chisquareall3reg}. The
$\chi^2$-values are somewhat large for small $N \le 100$, but once
$N\ge 200$ they are small for all families except for the connected
bipartite graphs, indicating good fits. For the connected bipartite
graphs, the $\chi^2$ values are small for $N$ large. This indicates
that perhaps the rate of convergence is slower for connected
bipartite graphs; we shall see additional differences in behavior
for these graphs below. Further, on average the $\chi^2$-values are
lowest for the $\beta=1$ case. While this suggests that the correct
model is a $\beta=1$ Tracy-Widom distribution, the data is not
conclusive.\\

\item \textbf{Percentage of eigenvalues to the left of the mean.} As
remarked, the four distributions, while close, differ in the
percentage of their mass to the left of their mean. By studying the
percentage of normalized eigenvalues in a sample less than the
sample mean, we see that the $\beta = 1$ distribution provides a
better fit to the observed results; however, with sample sizes of
1000 all four distributions provide good fits (see Table
\ref{table:massleft3INd}). We
therefore increased the number of graphs in the samples from 1000 to
100,000 for $N \in \{1002$, $2000$, $5002\}$ for the four families;
increasing the sample size by a factor of 100 gives us an additional
decimal digit of accuracy in measuring the percentages. See Table
\ref{table:massleft52test} for the results; this is \emph{the} most
important experiment in the paper, and shows that for the families
$\mathcal{CI}_{N,d}$, $\mathcal{SCI}_{N,d}$, and
$\mathcal{SCB}_{N,d}$ the $\beta=1$ Tracy-Widom distribution
provides a significant fit, but the other three distributions do
not. Thus we have found a statistic which is sensitive to very fine
differences between the four normalized distributions. However,
\emph{none} of the four candidate distributions provide a good fit
for the family $\mathcal{CB}_{N,d}$ for these values of $N$. For
this family the best fit is still with $\beta=1$, but the
$z$-statistics are high (between 3 and 4), which suggests that
either the distribution of eigenvalues for $d$-regular connected
bipartite graphs might not be given by a $\beta=1$ Tracy-Widom
distribution, or that the rate of convergence is slower; note our
$\chi^2$-tests suggests that the rate of convergence is indeed
slower for the connected bipartite family. In fact, upon increasing
$N$ to 10022 we obtain a good fit for connected bipartite graphs;
the $z$-statistic is about 2 for $\beta=1$, and almost 5 or larger
for the other three distributions. We shall see below that there are
other statistics where this family behaves differently than the
other three, strongly suggesting its rate of convergence
is slower.   \\

\item \textbf{Independence of $\lambda_\pm(G)$.} A graph is Ramanujan if $|\lambda_\pm(G)| \le 1$. For bipartite graphs it suffices to study $\lambda_+(G)$, as $\lambda_-(G) = -\lambda_+(G)$. For the non-bipartite families, however, we must investigate both. For our non-bipartite families we computed the sample correlation coefficient\footnote{The sample correlation coefficient $r_{xy}$ is $S_{xy} / \sqrt{S_{xx}S_{yy}}$, where $S_{uv} = \sum_{i=1}^n (u_i - \overline{u})(v_i-\overline{v})$ (with $\overline{u}$ the mean of the $u_i$'s). By Cauchy-Schwarz, $|r_{xy}| \le 1$. If the $x_i$ and $y_i$ are independent then $r_{xy}=0$, though the converse need not hold.} for $\lambda_+(G)$ and $\lambda_-(G)$ as $G$ varied through our random sample of 1000 graphs with $N$ vertices. For the $\mathcal{SCI}_{N,d}$ families we found the correlation coefficients were quite small; when $d=3$ they were in $[-.0355, 0.0827]$. For the $\mathcal{CI}_{N,d}$ the values were larger, but still small. When $d=3$ the correlation coefficients were in $[-0.0151, 0.2868]$, and all but two families with at least 5000 vertices had a correlation coefficient less than .1 in absolute value (and the values were generally decreasing with increasing $N$); see Figure \ref{fig:indepCISCI3} for the values. Thus the data suggests that $\lambda_\pm(G)$ are independent (for non-bipartite families). \\

%\ \\ \huge \textbf{I AM CURRENTLY GOING THROUGH AND FIXING THE %PAPER. I HAVE MADE IT TO THIS POINT.}\\ \normalsize

\item \textbf{Percentage of graphs that are Ramanujan.}
Except occasionally for the connected bipartite families, almost always
$s_\pm(\mathcal{F}_{N,d}) > m_\pm(\mathcal{F}_{N,d})$. Recall our
normalization of the eigenvalues from
\eqref{eq:norm2ndlargestevfam}: \be \widetilde{\gl_\pm}(G) \ = \
\frac{\gl_\pm(G) - 2\sqrt{d-1}+c_{\mu,N,d,\pm}
N^{m_\pm(\mathcal{F}_{N,d})}}{c_{\sigma,N,d,\pm} N^{s_\pm(\mathcal{F}_{N,d})}};
\ee Log-log plots of the differences between the sample means and
the predicted values, and standard deviations yield behavior that is
approximately linear as a function of $\log N$, supporting the
claimed normalization. Further, the exponents appear to be almost
constant in $N$, depending mostly only on $d$ (see Figure
\ref{fig:Mean3Reg}).
If this behavior holds as $N\to\infty$ then in the limit
approximately 52\% of the time we have $\gl_+(G) \le 2\sqrt{d-1}$ (and similarly about 52\% of the time $|\gl_-(G)| \le 2\sqrt{d-1}$. As $\gl_-(G) = -\gl_+(G)$ for bipartite graphs, this implies that about 52\% of the time bipartite graphs will be Ramanujan. Non-bipartite families behave differently. Assuming $\gl_+(G)$ and $\gl_-(G)$ are independent, the probability that both are at most $2\sqrt{d-1}$ in absolute value is about $27\%$ ($52\% \cdot 52\%$). See Figure
\ref{fig:raman3Reg} for plots of the percentages and Conjecture \ref{conj:32} for exact statements of these probabilities.
Unfortunately the rate of convergence is too slow
for us to see the conjectured limiting behavior.\\

\ei

Based on our results, we are led to the following conjecture.

\begin{conj}\label{conj:32} Let $\mathcal{F}_{N,d}$ be one of the
following families of $d$-regular graphs: $\mathcal{CI}_{N,d}$,
$\mathcal{SCI}_{N,d}$, or $\mathcal{SCB}_{N,d}$ (see Remark
\ref{rek:familiesdreggraphs} for definitions). The distribution of
$\gl_\pm(G)$, appropriately normalized as in
\eqref{eq:norm2ndlargestevfam}, converges as $N\to\infty$ to the
$\beta=1$ Tracy-Widom distribution (and not to a normalized
$\beta=2$ or $\beta=4$ Tracy-Widom distribution, or the standard normal distribution). For non-bipartite graphs, $\gl_+(G)$ and $\gl_-(G)$ are statistically independent. The normalization constants have $c_{\mu,N,d,\pm} < 0$ and $s_\pm(\mathcal{F}_{N,d}) > m_\pm(\mathcal{F}_{N,d})$, implying that in the limit as $N\to\infty$ approximately $52\%$ of the graphs in the bipartite families and 27\% otherwise are Ramanujan (i.e., $\gl(G) \le 2\sqrt{d-1}$); the actual percentage for the bipartite graphs is the percent of
mass in a $\beta=1$ Tracy-Widom distribution to the left of the mean
(to six digits it is 51.9652\%), and the square of this otherwise.
\end{conj}

\begin{rek} The evidence for the above conjecture is very strong for
three families. While the conjecture is likely to be true for
the connected bipartite graphs as well, different behavior is
observed for smaller $N$, though this may simply indicate a slower
rate of convergence. For example, when we
studied the percentage of eigenvalues to the left of the sample
mean, this was the only family where we did not obtain good fits to
the normalized $\beta=1$ Tracy-Widom distribution for $N \le 5002$,
though we did obtain good fits at $N=10022$ (see Table
\ref{table:massleft52test})
\end{rek}

\vfill

%%%%%%%%%%%%%%%%%%%%%%%%%%%%%%%%%%%%%%%%%%%%%%%%%%%%%%%%%%%%%%%%%%%%%%%%%%%%%%%%
%%%%%%%%%%%%%%%%%%%%%%%%%%%%%%%%%%%%%%%%%%%%%%%%%%%%%%%%%%%%%%%%%%%%%%%%%%%%%%%%

\pagebreak[4]
\newpage

\section{Results for $3$-Regular Graphs}\label{sec:3regstuff}

For $N \in \{$26, 32, 40, 50, 64, 80, 100, 126, 158, 200, 252, 316,
400, 502, 632, 796, 1002, 1262, 1588, 2000, 2516, 3168, 3990, 5022,
6324, 7962, 10022, 12618, 15886, 20000$\}$, we randomly chose 1000
$3$-regular graphs from the families $\mathcal{CI}_{N,3}$,
$\mathcal{SCI}_{N,3}$, $\mathcal{CB}_{N,3}$ and
$\mathcal{SCB}_{N,3}$. We analyzed the distributions of $\gl_\pm(G)$ for each sample using Matlab's eigs function\footnote{The Matlab code was originally written to investigate bipartite graphs. The symmetry of the eigenvalues allowed us to just look at the second largest eigenvalue; when we ran the code for non-bipartite graphs, we originally did not realize this had been hardwired. Thus there we were implicitly assuming $\lambda(G)=\lambda_+(G)$, which is frequently false for non-bipartite graphs. This error led us to initially conjecture 52\% of these graphs are Ramanujan in the limit, instead of the 27\% we discuss later.}, and investigated whether or not
it is well-modeled by the $\beta=1$ Tracy-Widom distribution.
Further, we calculated what percent of graphs were Ramanujan as well
as what percent of graphs had $|\gl_\pm(G)|$ less than
the sample mean; these statistics help elucidate the behavior as the
number of vertices tends to infinity.

%%%%%%%%%%%%%%%%%%%%%%%%%%%%%%%%%%%%%%%%%%%%%%%%%%%%%%%%%%%%%%%%%%%%%%%%%%%%%%%%%
\subsection{Distribution of $\gl_\pm(G)$}

\begin{figure}[ht]
\begin{center}
\caption{\label{fig:plotsofalldis3} Distribution of $\gl_+(G)$
for 1000 graphs randomly chosen from the ensemble
$\mathcal{CI}_{N,3}$ for various $N$. The vertical line is
$2\sqrt{2}$ and $N \in \{3990, 5022, 6324, 7962, 10022, 12618\}$. The curve with the lowest maximum value corresponds to $N = 3990$, and as $N$ increases the maximum value increases (so $N = 12618$ corresponds to the curve with greatest maximum value).} \scalebox{.4}{\includegraphics{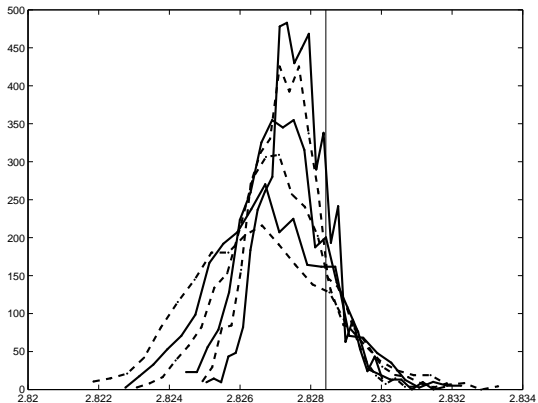}}
\end{center}\end{figure}

In Figure \ref{fig:plotsofalldis3} we plot the histogram
distribution of $\gl_+(G)$ for
$\mathcal{CI}_{N,3}$; the plots for the other families and for $\gl_-(G)$ are similar. This is a plot of
the actual eigenvalues. To determine whether or not the $\beta = 1$
Tracy-Widom distribution (or another value of $\beta$ or even a
normal distribution) gives a good fit to the data we rescale the
samples to have mean 0 and variance 1, and then compare the results
to scaled Tracy-Widom distributions (and the standard normal). In
Table \ref{table:chisquare3reg} we study the $\chi^2$-values for the
fits from the three Tracy-Widom distributions and the normal
distribution.

%\newpage

\newpage

\begin{table}[ht]
\caption{\label{table:chisquare3reg}Summary of $\chi^2$-values: each set is 1000 random $3$-regular graphs from
$\mathcal{CI}_{N,3}$ with $N \in \{$26, 32, 40, 50, 64, 80, 100, 126, 158, 200, 252, 316, 400, 502, 632, 796, 1002, 1262, 1588, 2000, 2516, 3168, 3990, 5022, 6324, 7962, 10022, 12618, 15886, 20000$\}$. The sample distribution in each set is
normalized to have mean 0 and variance 1, and is then compared to
normalized Tracy-Widom distributions ${\rm TW}_\beta^{{\rm norm}}$
($\beta \in \{1,2,4\}$, normalized to have mean 0 and variance 1)
and the standard normal ${\rm N}(0,1)$. There are 19 degrees of
freedom, and the critical values are 30.1435 (for $\alpha = .05$)
and 36.1908 (for $\alpha = .01$). }
\begin{center}
\begin{tabular}{|r|rrrr|}\hline
$N$  &   ${\rm TW}_1^{{\rm norm}}$  & ${\rm TW}_2^{{\rm norm}}$ &
${\rm TW}_4^{{\rm norm}}$ &  ${\rm N}(0,1)$\\ \hline
mean (all $N$) &   27.0	&	24.5	&	24.0	&	29.4	\\
median (all $N$) & 21.2	&	19.1	&	20.0	&	26.5	\\
%standard deviation (all)   &   42  &   18  &   180 &   7   \\
\hline
mean (last 10) &   21.7	&	22.2	&	23.7	&	35.0\\
median (last 10)& 21.2	&	20.9	&	22.4	&	35.4	\\
%standard deviation (last 10) &   8   &   5   &   37  &   8   \\
 \hline
%mean (last 5) &    22.1    &   22.7    &   22.2    &   41.6    \\
%median (last 5) & 21.7    &   19.8    &   19.3    &   27.3    \\
%\hline
 \end{tabular}
 \end{center}
\end{table}

\newpage

\begin{table}
\caption{\label{table:chisquareall3reg} Summary of $\chi^2$-values: each set is 1000 random $3$-regular graphs with $N$ vertices from our
families, with $N \in \{$26, 32, 40, 50, 64, 80, 100, 126, 158, 200, 252, 316, 400, 502, 632, 796, 1002, 1262, 1588, 2000, 2516, 3168, 3990, 5022, 6324, 7962, 10022, 12618, 15886, 20000$\}$. The sample distribution in each set is normalized to have
mean 0 and variance 1, and is then compared to the normalized
$\beta=1$ Tracy-Widom distributions. There are 19 degrees of
freedom, and the critical values are 30.1435 (for $\alpha = .05$)
and 36.1908 (for $\alpha = .01$).}
\begin{center}
\begin{tabular}{|r|rrrr|}\hline
$N$  &   $\mathcal{CI}_{N,3}$ &   $\mathcal{SCI}_{N,3}$  &
$\mathcal{CB}_{N,3}$ &  $\mathcal{SCB}_{N,3}$ \\ \hline
mean (all $N$)    &   27  &   19  &   78  &   19  \\
standard deviation (all $N$)   &   21  &   8  &   180 &   7   \\
\hline \hline
mean (last 10)               &   22  &   18  &   44  &   17  \\
standard deviation (last 10) &   11   &   6   &   37  &   8   \\
\hline \hline
mean (last 5)   &   23  &   18  &   32  &   14  \\
standard deviation (last 5)  &   13  &   8   &   23  &   1   \\
\hline
 \end{tabular}
 \end{center}
\end{table}

\newpage

As Table \ref{table:chisquare3reg} shows, the three normalized
Tracy-Widom distributions all give good fits, and even the standard
normal gives a reasonable fit.\footnote{While the data displayed above is for $\lambda_+(G)$, the $\chi^2$ values for $\lambda_-(G)$ and $\lambda(G)$ are comparable.} We divided the data into 20 bins and
calculated the $\chi^2$-values; with 19 degrees of freedom, the
$\alpha = .05$ threshold is 30.1435 and the $\alpha = .01$ threshold
is 36.1908.\footnote{We could use the (pessimistic) Bonferroni adjustments for multiple
comparisons (for ten comparisons these numbers become
38.5822 and 43.8201); we do not do this as the fits are already quite
good.} We investigate below another statistic which is
better able to distinguish the four candidate distributions. We note
that the normalized $\beta=1$ distribution gives good fits as
$N\to\infty$ for all the families, as indicated by Table
\ref{table:chisquareall3reg}. The fits are good for modest $N$ for
all families but the connected bipartite graphs; there the fit is
poor until $N$ is large. This indicates that the connected bipartite
graphs may have slower convergence properties than the other
families.

In Table \ref{table:paramsTW} we listed the mass to the left of the
mean for the Tracy-Widom distributions; it is 0.519652 for
$\beta=1$, 0.515016 for $\beta = 2$ and 0.511072 for $\beta=4$ (note
it is .5 for the standard normal). Thus looking at the mass to the
left of the sample mean provides a way to distinguish the four
candidate distributions; we present the results of these
computations for each set of 1000 graphs from $\mathcal{CI}_{N,3}$
in Table \ref{table:massleft3INd} (the other families behave
similarly). If $\theta_{{\rm obs}}$ is the observed percent of the
sample data (of size 1000) below the sample mean, then the
$z$-statistic \be z\ =\ (\theta_{{\rm obs}} - \theta_{{\rm pred}})
\Big/ \sqrt{\theta_{{\rm pred}} \cdot (1-\theta_{{\rm pred}}) /
1000} \ee measures whether or not the data supports that
$\theta_{{\rm pred}}$ is the percent below the mean.

%\newpage

\begin{table}[ht]
\caption{\label{table:massleft3INd} The mass to the left of the
sample mean for $\gl_+(G)$ from each set of 1000 $3$-regular graphs from
$\mathcal{CI}_{N,3}$ and the corresponding $z$-statistics comparing
that to the mass to the left of the mean of the three Tracy-Widom
distributions (0.519652 for $\beta = 1$, 0.515016 for $\beta=2$,
0.511072 for $\beta=4$) and the standard normal (.500). We use the
absolute value of the $z$-statistics for the means and medians. For
a two-sided $z$-test, the critical thresholds are 1.96 (for $\alpha
= .05$) and 2.575 (for $\alpha = .01$). For brevity we report only some of the values for $N \in \{$26, 32, 40, 50, 64, 80, 100, 126, 158, 200, 252, 316, 400, 502, 632, 796, 1002, 1262, 1588, 2000, 2516, 3168, 3990, 5022, 6324, 7962, 10022, 12618, 15886, 20000$\}$, but list the mean and medians for the last 5 and last 10 values of $N$.  }
\begin{center}
\begin{tabular}{|r|c||rrrr|}
  \hline
  $N$ & Observed mass &\ \ $z_{{\rm TW},1}$\ \ & $z_{{\rm TW},2}$ &\ \ $z_{{\rm TW},4}$\ \ &
  $z_{{\rm StdNorm}}$ \\ \hline
26	&	0.477	&	-2.700	&	-2.405	&	-2.155	&	-1.455	\\
100	&	0.522	&	0.149	&	0.442	&	0.691	&	1.391	\\
400	&	0.522	&	0.149	&	0.442	&	0.691	&	1.391	\\
1588	&	0.526	&	0.402	&	0.695	&	0.944	&	1.644	 \\
6324	&	0.524	&	0.275	&	0.568	&	0.818	&	1.518	 \\
20000	&	0.551	&	1.984	&	2.277	&	2.526	&	3.226	 \\
\hline \hline
mean (last 10)  &  0.519 & 0.861   &   0.873   &   0.960   &   1.341     \\
median (last 10) & 0.519   &   0.696   &   0.758   &   0.854   &   1.170   \\
\hline
mean (last 5)   &   0.514   &   1.186   &  1.126  &   1.076   &   1.138   \\
median (last 5) &   0.508   &   1.434   &   1.140   &   0.890   &   0.506   \\
\hline
\end{tabular}
\end{center}
\end{table}

%\newpage

While the data in Table \ref{table:massleft3INd} suggests that the
$\beta=1$ Tracy-Widom is the best fit, the other three distributions
provide good fits as well. As we expect the fit to
improve as $N$ increases, the last few rows of the table are the
most important. In 5 of the last 10 rows the smallest $z$-statistic
is with the $\beta=1$ Tracy-Widom distribution. Further, the average
of the absolute values of the $z$-values for the last 10 rows are 0.861 ($\beta=1$), 0.873
($\beta=2$), 0.960 ($\beta=4$) and 1.341 (for the standard normal),
again supporting the claim that the best fit is from the $\beta=1$
Tracy-Widom distribution.

In order to obtain more conclusive evidence as to which distribution
best models the second largest normalized eigenvalue, we considered
larger sample sizes (100,000 instead of 1000) for all four families;
see Table \ref{table:massleft52test} for the analysis. While there is a sizable increase in run-time (it took on the order of a few days to run the simulations for the three different values of $N$ for the
four families), we gain a decimal digit of precision in estimating
the percentages. This will allow us to statistically distinguish the
four candidate distributions.

%\pagebreak[4]

%\newpage 

\begin{table}[ht]
\caption{\label{table:massleft52test} The mass to the left of the
sample mean of $\gl_+(G)$ for each set of 100,000 $3$-regular graphs from our four
families ($\mathcal{CI}_{N,3}$,   $\mathcal{SCI}_{N,3}$,
$\mathcal{CB}_{N,3}$ and  $\mathcal{SCB}_{N,3}$), and the
corresponding $z$-statistics comparing that to the mass to the left
of the mean of the three Tracy-Widom distributions (0.519652 for
$\beta = 1$, 0.515016 for $\beta=2$, 0.511072 for $\beta=4$) and the
standard normal (.500). Discarded refers to the number of graphs
where Matlab's algorithm to determine the second largest eigenvalue
did not converge; this was never greater than 4 for any data set.
For a two-sided $z$-test, the critical thresholds are 1.96 (for
$\alpha = .05$) and 2.575 (for $\alpha = .01$).}
\begin{center}
\begin{tabular}{|c||rrrcc|}
  \hline
$\mathcal{CI}_{N,3}$   &   $z_{{\rm TW},1}$\ \ & $z_{{\rm TW},2}$ &\
\ $z_{{\rm TW},4}$\ \ &  $z_{{\rm StdNorm}}$  & Discarded   \\
\hline
1002	&	1.2773	&	4.2103	&	6.7044	&	13.7053	& 0 \\
2000	&	0.9671	&	3.9002	&	6.3944	&	13.3954	& 0 \\
5022	&	0.3152	&	3.2485	&	5.7428	&	12.744 & 0	\\
\  &\ &\ &\ & \  \\ \hline
 $\mathcal{SCI}_{N,3}$   &   $z_{{\rm TW},1}$\ \ &
$z_{{\rm TW},2}$ &\ \ $z_{{\rm TW},4}$\ \ &  $z_{{\rm StdNorm}}$  &
Discarded
\\ \hline
1002	&	-0.7481	&	2.1855	&	4.6801	&	11.6815 & 0	\\
2000	&	-0.5899	&	2.3437	&	4.8382	&	11.8396 & 0	\\
5022	&	-1.0456	&	1.8881	&	4.3827	&	11.3842 & 0	\\
\  &\ &\ &\ & \  \\ \hline
 $\mathcal{CB}_{N,3}$   &   $z_{{\rm
TW},1}$\ \ & $z_{{\rm TW},2}$ &\ \ $z_{{\rm TW},4}$\ \ &  $z_{{\rm
StdNorm}}$  & Discarded   \\ \hline
1002    &   3.151   &   6.083   &   8.577   &   15.577  &   0   \\
2000    &   3.787   &   6.719   &   9.213   &   16.213  &   1   \\
5022    &   3.563   &   6.495   &   8.989   &   15.989  &   4   \\
10022\ \  &   2.049   &   4.982   &   7.476   &   14.477  &   0   \\
\  &\ &\ &\ & \  \\ \hline
 $\mathcal{SCB}_{N,3}$   &   $z_{{\rm
TW},1}$\ \ & $z_{{\rm TW},2}$ &\ \ $z_{{\rm TW},4}$\ \ &  $z_{{\rm
StdNorm}}$  & Discarded   \\ \hline
1002    &   -1.963  &   0.971   &   3.465   &   10.467  &   0   \\
2000    &   -0.767  &   2.167   &   4.661   &   11.663  &   2   \\
5022    &   -0.064  &   2.869   &   5.364   &   12.365  &   4   \\
\hline
\end{tabular}
\end{center}
\end{table}

%\newpage

\emph{This is the most important test in the paper.} The results are
striking, and strongly support that only the $\beta = 1$ Tracy-Widom
distribution models $\gl_\pm(G)$ (the results for $\lambda_-(G)$ were similar to those for $\lambda_+(G)$). Except for
$\mathcal{SCB}_{1002,3}$, for each of the families and each $N$ the
$z$-statistic increases in absolute value as we move from $\beta=1$
to $\beta=2$ to $\beta=4$ to the standard normal. Further, the
$z$-values indicate excellent fits with the $\beta=1$ distribution
for all $N$ and all families \emph{except} the $3$-regular connected
bipartite graphs; no other value of $\beta$ or the standard normal
give as good of a fit. In fact, the other fits are often terrible.
The $\beta=4$ and standard normal typically have $z$-values greater
than $4$; the $\beta=2$ gives a better fit, but significantly worse
than $\beta = 1$. 

Thus, except for $3$-regular connected bipartite graphs, the data is
consistent only with a $\beta =1$ Tracy-Widom distribution. In the
next subsections we shall study the sample means, standard
deviations, and percent of graphs in a family that are Ramanujan. We
shall see that the $3$-regular connected bipartite graphs
consistently behave differently than the other three families (see
in particular Figure \ref{fig:raman3Reg}).

\pagebreak[4]
\newpage

%%%%%%%%%%%%%%%%%%%%%%%%%%%%%%%%%%%%%%%%%%%%%%%%%%%%%%%%%%%%%%%%%%%%%%%%%%%%%%%%%
\subsection{Means and Standard Deviations}

\begin{figure}
\begin{center}
\caption{\label{fig:Mean3Reg} Sample means of $\lambda_+(G)$: each set is 1000 random
$3$-regular graphs with $N$ vertices, chosen according to the
specified construction. The first plot is the mean versus the number
of vertices; the second plot is a log-log plot of the mean and the
number of vertices. $\mathcal{CI}_{N,3}$ are stars,
$\mathcal{SCI}_{N,3}$ are triangles, $\mathcal{CB}_{N,3}$ are diamonds,
$\mathcal{SCB}_{N,3}$ are boxes; the dashed line is $2 \sqrt{2}
\approx 2.8284$.} \scalebox{.6}{\includegraphics{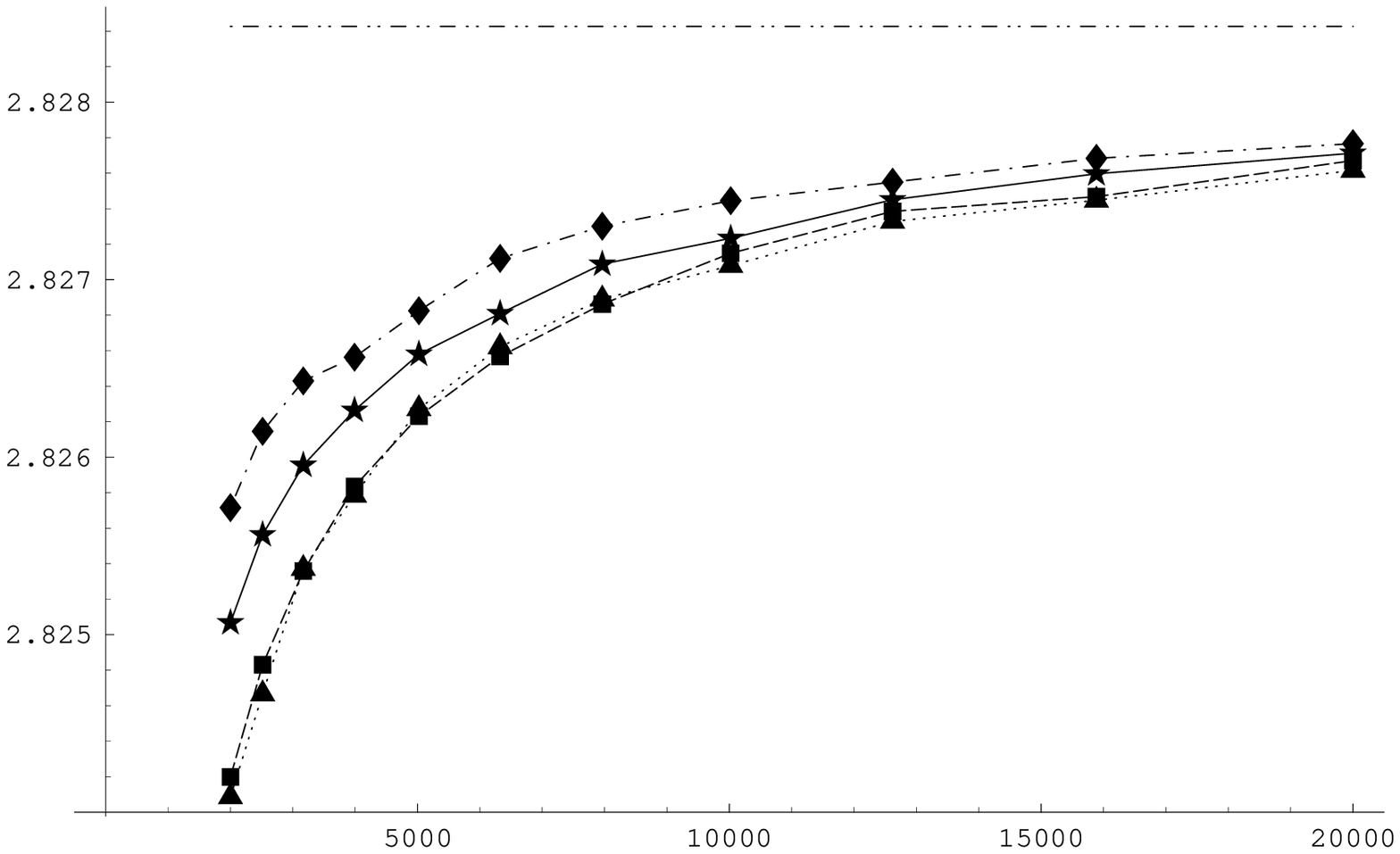}}
\scalebox{.6}{\includegraphics{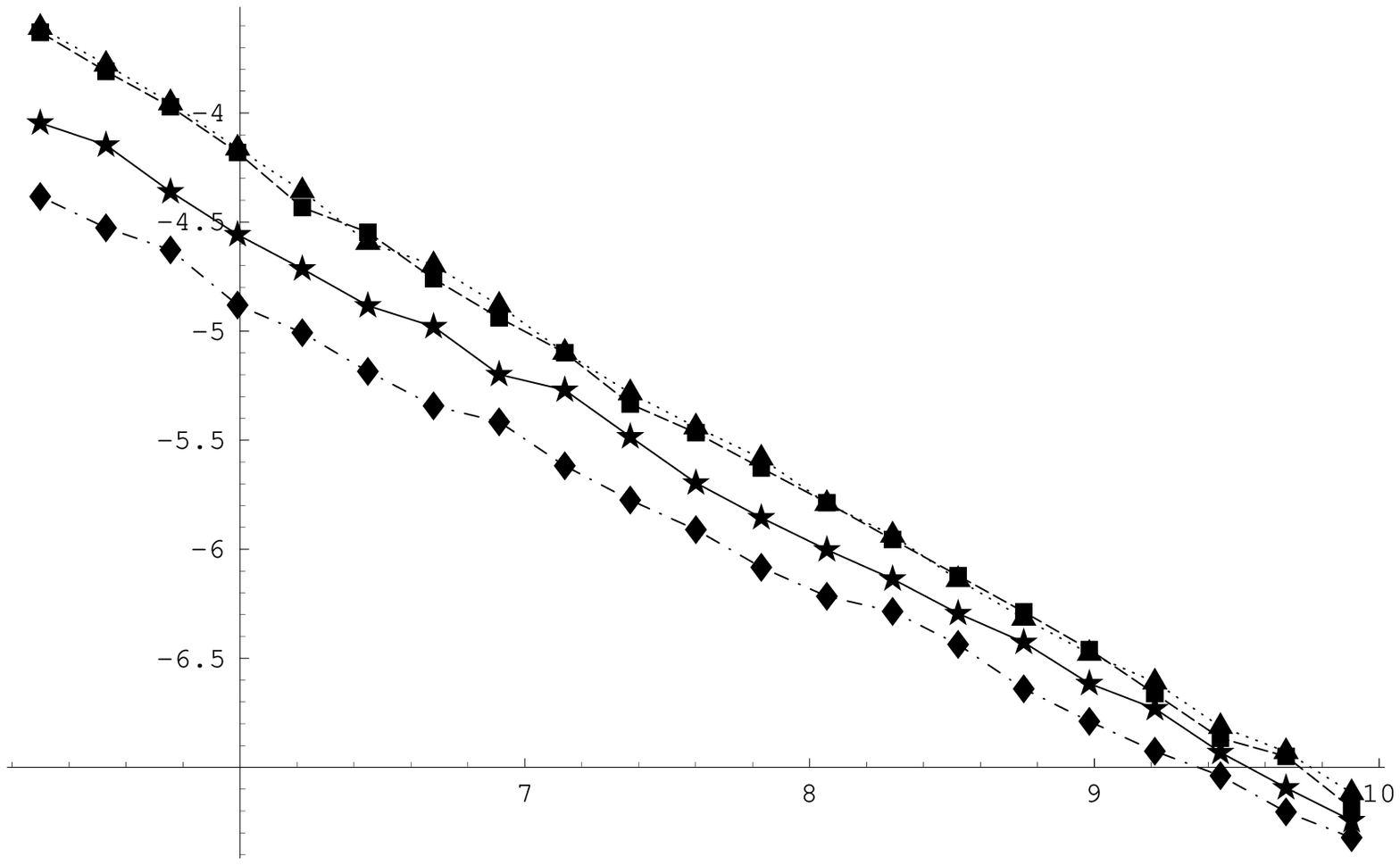}}
\end{center}\end{figure}

In Figure \ref{fig:Mean3Reg} we plot the sample means of sets of
1000 $3$-regular graphs chosen randomly from $\mathcal{CI}_{N,3}$
(connected perfect matchings), $\mathcal{SCI}_{N,3}$ (simple
connected perfect matchings), $\mathcal{CB}_{N,3}$  (connected
bipartite) and  $\mathcal{SCB}_{N,3}$ (simple connected bipartite) against the number of vertices.

Because of analogies with similar systems whose largest eigenvalue
satisfies a Tracy-Widom distribution, we expect the normalization
factor for the second largest eigenvalue to be similar to that in
\eqref{eq:normlargevTW1gauss}. As we do not expect that the factors
will still be $N^{1/2}$ and $N^{1/6}$, we consider the general
normalization given in \eqref{eq:norm2ndlargestevfam}; for a
$3$-regular graph in one of our families we study \be
\widetilde{\gl_\pm}(G) \ = \ \frac{|\gl_\pm(G)| - 2\sqrt{2}+c_{\mu,N,3,\pm}
N^{m_\pm(\mathcal{F}_{N,3})}}{c_{\sigma,N,3,\pm} N^{s_\pm(\mathcal{F}_{N,3})}}.
\ee

\begin{rek}\label{rek:whysmdiffmatters}
The most important parameters are the exponents
$m_\pm(\mathcal{F}_{N,3})$ and $s_\pm(\mathcal{F}_{N,3})$; previous work
\cite{Fr2} (and our investigations) suggest that $c_{\mu,N,3,\pm} < 0$.
Let us assume that, in the limit as the number of vertices tends to
infinity, the distributions of $|\lambda_\pm(G)|$ converge to the $\beta = 1$ Tracy-Widom distribution and
that $c_{\mu,N,3,\pm} < 0$. If $s_\pm(\mathcal{F}_{N,3}) >
m_\pm(\mathcal{F}_{N,3})$ then in the limit we expect about 52\% of the
graphs to have $\gl_+(G) \le 2\sqrt{2}$ (and similarly for $|\gl_-(G)|$), as
this is the mass of the $\beta=1$ Tracy-Widom distribution to the
left of the mean. To see why this is true, note that if
$\mu^{}_{\mathcal{F}_{N,3,+}}$ and $\sigma^{}_{\mathcal{F}_{N,3,+}}$ are
the mean and standard deviation of the data set of $\lambda_+(G)$
for all $G\in\mathcal{F}_{N,3}$, then
$\mu^{}_{\mathcal{F}_{N,3,+}}\approx2\sqrt{2} -
c_{\mu,N,d,+}N^{m_+(\mathcal{F}_{N,3})}$ and
$\sigma^{}_{\mathcal{F}_{N,3,+}}\approx
c_{\sigma,N,3,+}N^{s_+(\mathcal{F}_{N,3})}$, so
 \be 2\sqrt{2}\ \approx\ \mu_{\mathcal{F}_{N,3,+}} +
\frac{c_{\mu,N,3,+}}{c_{\sigma,N,3,+}} \cdot
N^{m_+(\mathcal{F}_{N,3})-s_+(\mathcal{F}_{N,3})} \cdot
\sigma_{\mathcal{F}_{N,3,+}}.\ee Thus the Ramanujan threshold,
$2\sqrt{2}$, will fall approximately $
\frac{c_{\mu,N,3,+}}{c_{\sigma,N,3,+}}N^{m_+(\mathcal{F}_{N,3})-s_+(\mathcal{F}_{N,3})}$
standard deviations away from the mean. In the limit as $N$ goes to
infinity we see that the threshold falls zero standard deviations to
the right of the mean if
$m_+(\mathcal{F}_{N,3})<s_+(\mathcal{F}_{N,3})$, but infinitely many if
$m_+(\mathcal{F}_{N,3})>s_+(\mathcal{F}_{N,3})$.
\end{rek}

%This is because, in this case, a very small multiple of the standard
%deviation suffices to bring us above $2\sqrt{2}$. If, however,
%$s(\mathcal{F}_{N,3}) < m(\mathcal{F}_{N,3})$ then in the limit we
%expect all graphs to be Ramanujan; this is because we are many
%standard deviations below the $2\sqrt{2}$.

\begin{table}[ht]
\caption{\label{table:bestfitmeansstd3Reg} The graph sizes are
chosen from $\{$26, 32, 40, 50, 64, 80, 100, 126, 158, 200, 252,
316, 400, 502, 632, 796, 1002, 1262, 1588, 2000, 2516, 3168, 3990,
5022, 6324, 7962, 10022, 12618, 15886, 20000$\}$. The first four
columns are the best-fit values of $m(\mathcal{F}_{N,3})$; the last
four columns are the best fit values of $s(\mathcal{F}_{N,3})$. Bold
entries are those where $s(\mathcal{F}_{N,3}) <
m(\mathcal{F}_{N,3})$; all other entries are where
$s(\mathcal{F}_{N,3}) > m(\mathcal{F}_{N,3})$. }
\begin{center}
\begin{tabular}{|r|cccc||cccc|}
  \hline
$N$  &   $\mathcal{CI}_{N,3}$ &   $\mathcal{SCI}_{N,3}$  &
$\mathcal{CB}_{N,3}$ &  $\mathcal{SCB}_{N,3}$  &
$\mathcal{CI}_{N,3}$ &   $\mathcal{SCI}_{N,3}$  &
$\mathcal{CB}_{N,3}$ &  $\mathcal{SCB}_{N,3}$    \\
\hline
$\{26,\dots,20000\}$    &   -0.792  &   -0.830  &   -0.723  &   -0.833  &   -0.718  &   -0.722  &   -0.709  &   -0.729  \\
$\{80,\dots,20000\}$    &   -0.756  &   -0.790  &   -0.671  &   -0.789  &   -0.701  &   -0.700  &   \textbf{-0.697}  &   -0.706  \\
$\{252,\dots,20000\}$   &   -0.727  &   -0.761  &   -0.638  &   -0.761  &   -0.695  &   -0.688  &   \textbf{-0.688}  &   -0.696  \\
$\{26,\dots,64\}$       &   -1.045  &   -1.097  &   -1.065  &   -1.151  &   -0.863  &   -0.906  &   -0.794  &   -0.957  \\
$\{80,\dots,200\}$      &   -0.887  &   -0.982  &   -0.982  &   -0.968  &   -0.769  &   -0.717  &   -0.719  &   -0.750  \\
$\{232,\dots,632\}$     &   -0.801  &   -0.885  &   -0.737  &   -0.842  &   -0.688  &   -0.713  &   -0.714  &   -0.734  \\
$\{796,\dots,2000\}$    &   -0.771  &   -0.819  &   -0.649  &   -0.785  &   -0.606  &   -0.719  &   \textbf{-0.705}  &   -0.763  \\
$\{2516,\dots,6324\}$   &   -0.745  &   -0.788  &   -0.579  &   -0.718  &   -0.714  &   -0.671  &   \textbf{-0.770}  &   -0.688  \\
$\{7962,\dots,20000\}$  &   -0.719  &   -0.692  &   -0.584  &   -0.757  &   -0.592  & \textbf{-0.707}  &   \textbf{-0.671}  &   -0.648  \\
\hline\end{tabular}
\end{center}
\end{table}

We record (some of) the best fit exponents in Table
\ref{table:bestfitmeansstd3Reg}; the remaining values are similar. To simplify the calculations, we
changed variables and did a log-log plot. Several trends can be seen
from the best fit exponents in Table
\ref{table:bestfitmeansstd3Reg}. Most of the time,
$s_\pm(\mathcal{F}_{N,3}) > m_\pm(\mathcal{F}_{N,3})$, which indicates that
it is more likely in the limit that 52\% (and not all) of the bipartite graphs
are Ramanujan (and about 27\% of the non-bipartite). Except for $\mathcal{CB}_{N,3}$ (connected bipartite
graphs), only once is $s_+(\mathcal{F}_{N,3}) < m_+(\mathcal{F}_{N,3})$;
for $\mathcal{CB}_{N,3}$ we have $s_+(\mathcal{F}_{N,3}) <
m_+(\mathcal{F}_{N,3})$ approximately half of the time. Further, the
best fit exponents $s_+(\mathcal{F}_{N,3})$ and $m_+(\mathcal{F}_{N,3})$
are mostly monotonically increasing with increasing $N$ (remember
all exponents are negative), and $c_{\mu,N,3,+}$ and $c_{\sigma,N,3,+}$
do not seem to get too large or small (these are the least important
of the parameters, and are dwarfed by the exponents). This suggests
that either the relationship is more complicated than we have
modeled, or $N$ is not large enough to see the limiting behavior.
While our largest $N$ is 20000, $\log(20000)$ is only about 10. Thus
we may not have gone far enough to see the true behavior. If the
correct parameter is $\log N$, it is unlikely that larger
simulations will help.

\begin{figure}
\begin{center}
\caption{\label{fig:loglogplotmean3regCIN3} Dependence of the
logarithm of the mean of $\gl_+(G)$ on $\log\left(-c_{\mu,N,3,+}
N^{m_+(\mathcal{CI}_{N,3})}\right)$ on $N$, showing the best fit lines using all 30 values of $N$ as well as just the last 10 values.}
\scalebox{.6}{\includegraphics{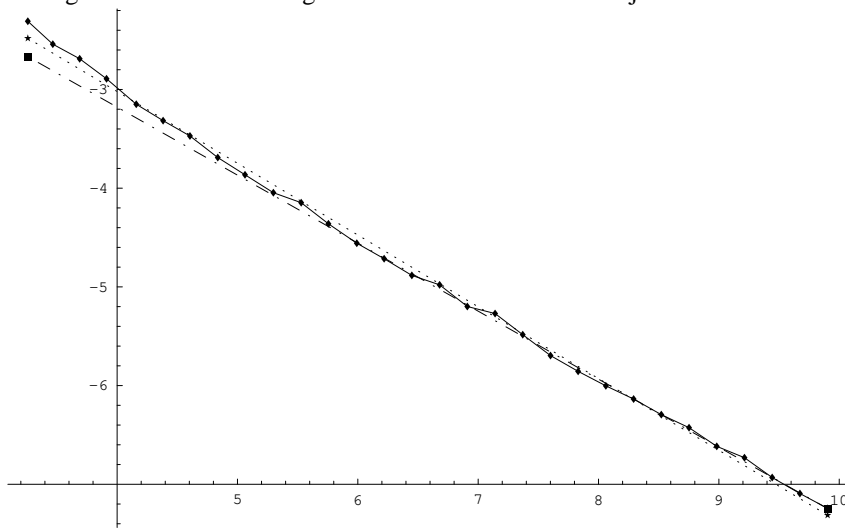}}
\end{center}\end{figure}

In Figure \ref{fig:loglogplotmean3regCIN3} we plot the
$N$-dependence of the logarithm of the difference of the mean from
$2\sqrt{2}$ versus the logarithm of $-c_{\mu,N,3,+}
N^{m_+(\mathcal{CI}_{N,3})}$, as well as the best fit lines obtained
by using all of the data and just the last 10 data points. As the
plot shows, the slope of the best fit line (the key parameter for
our investigations) noticeably changes in the region we investigate,
suggesting that either we have not gone high enough to see the limiting,
asymptotic behavior or that it is not precisely linear.

%%%%%%%%%%%%%%%%%%%%%%%%%%%%%%%%%%%%%%%%%%%%%%%%%%%%%%%%%%%%%%%%%%%%%%%%%%%%%%%%%
\subsection{Independence of $\gl_\pm(G)$ in non-bipartite families}

In determining what percentage of graphs in a non-bipartite family is Ramanujan, it is important to know whether or not $\gl_+(G)$ and $\gl_-(G)$ are statistically independent as $G$ varies in a family. For example, if they are perfectly correlated the percentage could be 100\%, while if they are perfectly anti-correlated it could be 0\%.

\begin{figure}
\begin{center}
\caption{\label{fig:indepCISCI3} Sample correlation coefficients of $\gl_\pm(G)$: each set is 1000
random $3$-regular graphs with $N$ vertices, chosen according to the
specified construction. We plot the sample correlation coefficient versus the
logarithm of the number of vertices. $\mathcal{CI}_{N,3}$ are stars and $\mathcal{SCI}_{N,3}$ are triangles.}
\scalebox{.6}{\includegraphics{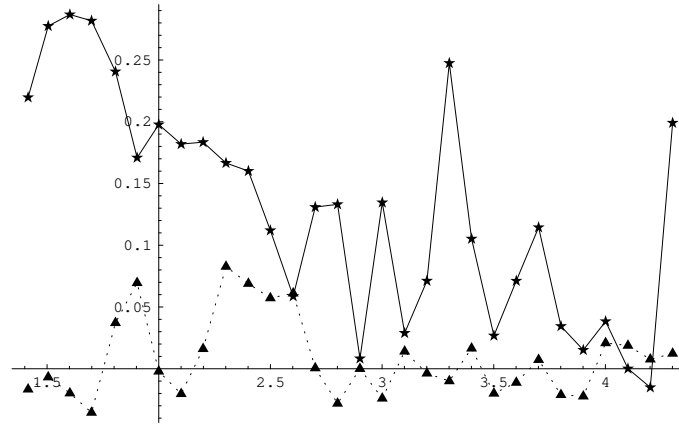}}
\end{center}\end{figure}

In Figure \ref{fig:indepCISCI3} we plot the sample correlation coefficient for $\gl_\pm(G)$ for the non-bipartite families. For $\mathcal{CI}_{N,3}$ the values are generally positive and decreasing with increasing $N$; for $\mathcal{SCI}_{N,3}$ the data appears uncorrelated, with very small coefficients oscillating about zero. As another test we compared the product of the observed probabilities that $\gl_+(G) < 2\sqrt{2}$ and $|\gl_-(G)| < 2\sqrt{2}$ to the observed probability that $\gl(G) < 2\sqrt{2}$; these values were virtually identical, which is what we would expect if $\gl_\pm(G)$ are statistically independent.

%%%%%%%%%%%%%%%%%%%%%%%%%%%%%%%%%%%%%%%%%%%%%%%%%%%%%%%%%%%%%%%%%%%%%%%%%%%%%%%%%
\subsection{Percentage of graphs that are Ramanujan}

\begin{figure}
\begin{center}
\caption{\label{fig:raman3Reg} Percent Ramanujan: each set is 1000
random $3$-regular graphs with $N$ vertices, chosen according to the
specified construction. The first plot is the percent versus the
number of vertices; the second plot is the percent versus the
logarithm of the number of vertices. $\mathcal{CI}_{N,3}$ are stars,
$\mathcal{SCI}_{N,3}$ are diamonds, $\mathcal{CB}_{N,3}$ are triangles,
$\mathcal{SCB}_{N,3}$ are boxes. }
\scalebox{.45}{\includegraphics{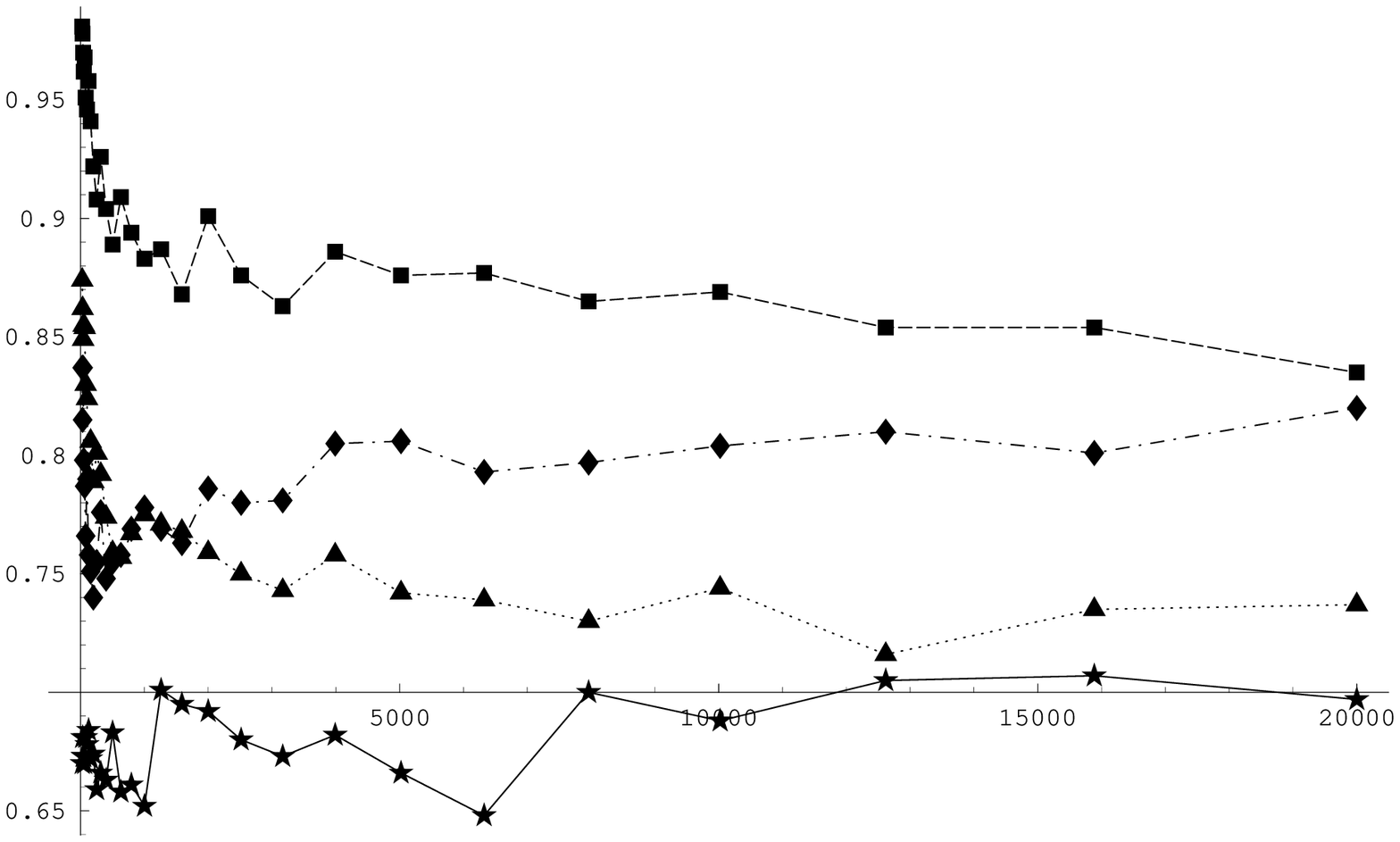}}\ \\
\scalebox{.45}{\includegraphics{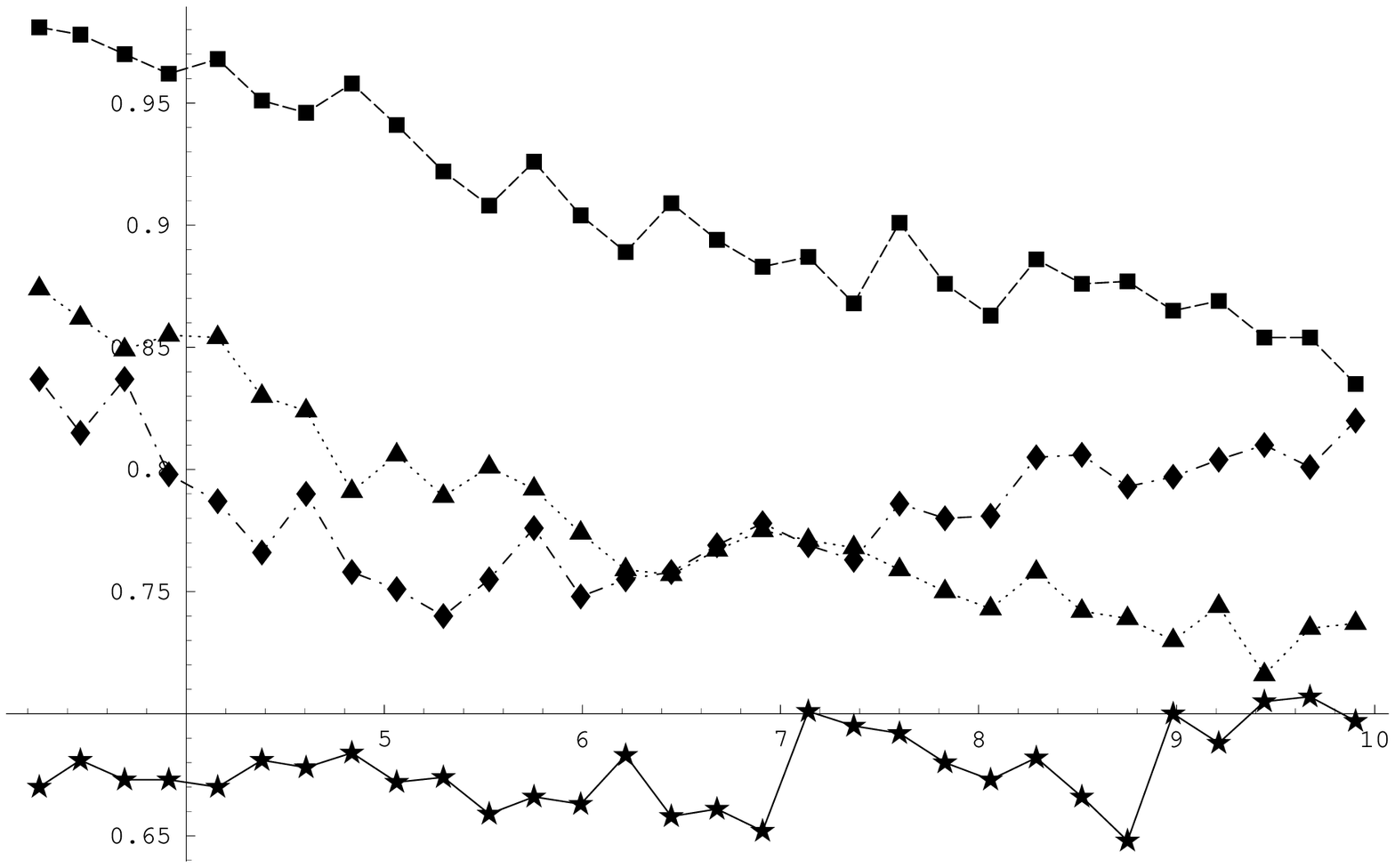}}
\end{center}\end{figure}

In Figure \ref{fig:raman3Reg} we plot the percentage of graphs in each
sample of 1000 from the four families that are Ramanujan (the first plot is
the percentage against the number of vertices, the second is the
percentage against the logarithm of the number of vertices).
The
most interesting observation is that, for the most part, the
probability that a random graph from the bipartite families is Ramanujan is
decreasing as $N$ increases, while the probability that a random
graph from the non-bipartite families is Ramanujan is oscillating in the range.
%As we saw in Table \ref{table:bestfitmeansstd3Reg} and
%Figure \ref{fig:loglogplotmean3regCIN3}, except for
%$\mathcal{CB}_{N,3}$ we had $s(\mathcal{F}_{N,3}) >
%m(\mathcal{F}_{N,3})$; for $\mathcal{CB}_{N,3}$ we only had
%$s(\mathcal{F}_{N,3}) > m(\mathcal{F}_{N,3})$ about half of the
%time.

\newpage

\newpage

\pagebreak[4]

%%%%%%%%%%%%%%%%%%%%%%%%%%%%%%%%%%%%%%%%%%%%%%%%%%%%%%%%%%%%%%%%%%%%%%%%%%%%%%%%
%%%%%%%%%%%%%%%%%%%%%%%%%%%%%%%%%%%%%%%%%%%%%%%%%%%%%%%%%%%%%%%%%%%%%%%%%%%%%%%%

\section*{Acknowledgement}

We thank Alex Barnett, Jon Bober, Peter Sarnak, Brad Weir, and the participants of the AMS Special Session on Expanders and Ramanujan Graphs: Construction and Applications at the 2008 Joint Meetings (organized by Michael T. Krebs,
Anthony M. Shaheen and
Audrey A. Terras)  for many enlightening discussions, Craig Tracy for sharing Mathematica code to compute the Tracy-Widom distributions, and the Information Technology Managers at the Mathematics Departments at Princeton, the Courant Institute and
Brown University for help in getting all the programs to run
compatibly. The first named author was partly supported by NSF grant
DMS0600848.

%%%%%%%%%%%%%%%%%%%%%%%%%%%%%%%%%%%%%%%%%%%%%%%%%%%%%%%%%%%%%%%%%%%%%%%%%%%%%%%%
%%%%%%%%%%%%%%%%%%%%%%%%%%%%%%%%%%%%%%%%%%%%%%%%%%%%%%%%%%%%%%%%%%%%%%%%%%%%%%%%

\bigskip \ \\ \bigskip


\begin{thebibliography}{999999}

\bibitem[Alh]{Alh}
Y. Alhassid, \emph{The statistical theory of quantum dots}, Rev.
Mod. Phys. \textbf{72} (2000), no. 4, 895--968.

\bibitem[Al]{Al}
N. Alon, \emph{Eigenvalues and expanders}, Combinatorica \textbf{6}
(1986), no. 2, 83--96.

\bibitem[AM]{AM}
N. Alon and V. Milman, \emph{$\lambda_1$, isoperimetric inequalities
for graphs, and superconcentrators},  J. Combin. Theory Ser. B
\textbf{38} (1985),  no. 1, 73--88.

\bibitem[BDJ]{BDJ}
J. Baik, P. Deift and K. Johansson, \emph{On the distribution of the
length of the longest increasing subsequence of random
permutations}, J. Amer. Math. Soc. \textbf{12} (1999), 1119--1178.

\bibitem[BR1]{BR1}
J. Baik and E. M. Rains, \emph{The asymptotics of monotone
subsequences of involutions}, Duke Math. J. \textbf{109} (2001),
205--281.

\bibitem[BR2]{BR2}
J. Baik and E. M. Rains, \emph{Symmetrized random permutations, in
Random Matrix Models and their Applications}, eds. P. Bleher and A.
Its, Math. Sci. Res. Inst. Publications \textbf{40}, Cambridge Univ.
Press, 2001, 1--19.

\bibitem[BR3]{BR3}
J. Baik and E. M. Rains, \emph{Limiting distributions for a
polynuclear growth model}, J. Stat. Phys. \textbf{100} (2000),
523--541.

\bibitem[Ba]{Ba}
Yu. Baryshnikov, \emph{GUEs and queues}, Probab. Th. Rel. Fields
\textbf{119} (2001), 256--274.

\bibitem[Bien]{Bien}
F. Bien, \emph{Constructions of telephone networks by group
representations}, Notices of the AMS \textbf{36} (1989), no. 1,
5--22.

\bibitem[Bol]{Bol}
B. Bollob\'{a}s, \emph{Random Graphs}, Cambridge Studies in
Advanced Mathematics, Cambridge University Press, 2001.

\bibitem[BOO]{BOO}
A. Borodin, A. Okounkov and G. Olshanski, \emph{Asymptotics of
Plancherel measures for symmetric groups}, J. Amer. Math. Soc.
\textbf{13} (2000), 481--515.

\bibitem[Chiu]{Chiu}
\newblock P. Chiu, \emph{Cubic Ramanujan graphs}, Combinatorica
\textbf{12} (1992), no. 3, 275--285.

\bibitem[DSV]{DSV}
G. Davidoff, P. Sarnak, and A. Valette, \emph{Elementary Number
Theory, Group Theory, and Ramanujan Graphs}, London Mathematical
Society, Student Texts \textbf{55}, Cambridge University Press,
2003.

\bibitem[Do]{Do}
J. Dodziuk, \emph{Difference equations, isoperimetric inequality and
transience of certain random walks},  Trans. Amer. Math. Soc.
\textbf{284} (1984),  no. 2, 787--794.

\bibitem[Fr1]{Fr1}
J. Friedman, \emph{Some geometric aspects of graphs and their eigenfunctions},   Duke Math. J. \textbf{69} (1993), no. 3, 487--525.


\bibitem[Fr2]{Fr2}
J. Friedman, \emph{A proof of Alon's second eigenvalue conjecture},
Proceedings of the Thirty-Fifth Annual ACM Symposium on Theory of
Computing,  720--724 (electronic), ACM, New York, 2003.

\bibitem[Gau]{Gau}
M. Gaudin, \emph{Sur la loi limite de l'espacement des valeurs
propres d'une matrice al\'{e}atoire}, Nucl. Phys. \textbf{25} (1961)
447--458.

\bibitem[GTW]{GTW}
\newblock J. Gravner, C. A. Tracy and H. Widom, \emph{Limit theorems for height
fluctuations in a class of discrete space and time growth models},
J. Stat. Phys. \textbf{102} (2001), 1085--1132.

\bibitem[GILVZ]{GILVZ}
O. Goldreich, R. Impagliazzo, L. Levin, R. Venkatesan, and D.
Zuckerman, \emph{Security preserving amplification of hardness}. In
31st Annual Symposium on Foundations of Computer Science, Vol. I, II
(St. Louis, MO, 1990), 318--326, IEEE Comput. Soc. Press, Los
Alamitos, CA, 1990.

\bibitem[HLW]{HLW}
S. Hoory, N. Linial and A. Wigderson, \emph{Expander graphs and their applications},   Bull. Amer. Math. Soc. \textbf{43} (2006), 439--561.

\bibitem[JMRR]{JMRR}
D. Jakobson, S. D. Miller, I. Rivin, and Z. Rudnick,
\emph{Eigenvalue spacings for regular graphs}. Pages 317--327 in
\emph{Emerging Applications of Number Theory (Minneapolis, 1996)},
The IMA Volumes in Mathematics and its Applications, Vol. 109,
Springer, New York, 1999.

%\bibitem[JLR]{JLR}
%S. Janson, T. Luczak, and A. Rucinski, \emph{Random Graphs}, John
%Wiley \& Sons, Inc., New York 2000.

\bibitem[Jo1]{Jo1}
K. Johansson, \emph{Discrete orthogonal polynomial ensembles and the
Plancherel measure}, Ann. Math. \textbf{153} (2001), 259--296.

\bibitem[Jo2]{Jo2}
K. Johansson, \emph{Non-intersecting paths, random tilings and
random matrices}, Probab. Th. Rel. Fields \textbf{123} (2002),
225--280.

\bibitem[Jo3]{Jo3}
K. Johansson, \emph{Toeplitz determinants, random growth and
determinantal processes}, ICM Vol. III (2002), 53--62.

\bibitem[KeSn]{KeSn}
\newblock J. P. Keating and N. C. Snaith, \emph{Random matrix theory
and $\zeta(1/2+it)$},  Comm. Math. Phys.  \textbf{214} (2000),  no.
1, 57--89.

\bibitem[LPS]{LPS}
A. Lubotzky, R. Phillips, and P. Sarnak, \emph{Ramanujan graphs},
Combinatorica  \textbf{8} (1988),  no. 3, 261--277.

\bibitem[Mar]{Mar}
G. A. Margulis, \emph{Explicit group-theoretic constructions of
combinatorial schemes and their applications in the construction of
expanders and concentrators (Russian)}, Problemy Peredachi
Informatsii \textbf{24} (1988), no. 1, 51--60;  translation in
Problems Inform. Transmission  \textbf{24} (1988),  no. 1, 39--46.

\bibitem[McK]{McK}
B. McKay, \emph{The expected eigenvalue distribution of a large
regular graph}, Linear Algebra Appl. \textbf{40} (1981), 203--216.

\bibitem[Meh1]{Meh1}
M. Mehta, \emph{On the statistical properties of level spacings in
nuclear spectra}, Nucl. Phys. \textbf{18} (1960), 395--419.

\bibitem[Meh2]{Meh2}
M. Mehta, \emph{Random Matrices}, 2nd edition, Academic Press,
Boston, 1991.

\bibitem[Mor]{Mor}
M. Morgenstern, \emph{Existence and explicit constructions of $q+1$
regular Ramanujan graphs for every prime power $q$},  J. Combin.
Theory Ser. B  \textbf{62}  (1994),  no. 1, 44--62.

\bibitem[Mur]{Mur}
M. Ram Murty, \emph{Ramanujan graphs}, J. Ramanujan Math. Soc.
\textbf{18} (2003), no. 1, 33--52.

\bibitem[Od1]{Od1}
\newblock A. Odlyzko, \emph{On the distribution of spacings
between zeros of the zeta function}, Math. Comp. \textbf{48} (1987),
no. 177, 273--308.

\bibitem[Od2]{Od2}
\newblock A. Odlyzko, \emph{The $10^{22}$-nd zero of the Riemann zeta function}, Proc.
Conference on Dynamical, Spectral and Arithmetic Zeta-Functions, M.
van Frankenhuysen and M. L. Lapidus, eds., Amer. Math. Soc.,
Contemporary Math. series, 2001,
http://www.research.att.com/$\sim$amo/doc/zeta.html

\bibitem[Pi]{Pi}
Pippenger, \emph{Super concentrators}, SIAM Journal Comp. \textbf{6}
(1977), 298--304.

\bibitem[PS1]{PS1}
M. Pr\"ahofer and H. Spohn, \emph{Statistical self-similarity of
one-dimensional growth processes}, Physica A \textbf{279} (2000),
342--352.

\bibitem[PS2]{PS2}
M. Pr\"ahofer and H. Spohn, \emph{Universal distributions for growth
processes in $1+1$ dimensions and random matrices}, Phys. Rev.
Letts. \textbf{84} (2000), 4882--4885.

\bibitem[Sar1]{Sar1}
P. Sarnak  \emph{Some applications of modular forms}, Cambridge
Trusts in Mathemetics, Vol. 99, Cambridge University Press,
Cambridge, 1990.

\bibitem[Sar2]{Sar2}
P. Sarnak  \emph{What is an Expander?}, Notices of the AMS \textbf{51} (2004), no. 7, 762--763.

\bibitem[So]{So}
A. Soshnikov, \emph{A note on universality of the distribution of
the largest eigenvalue in certain classes of sample covariance
matrices}, preprint (arXiv: math.PR/0104113).

\bibitem[SS]{SS}
M. Sipser and D. A. Spielman, \emph{Expander codes}, IEEE Trans.
Inform. Theory  \textbf{42}  (1996),  no. 6, part 1, 1710--1722.

\bibitem[TW1]{TW1}
C. A. Tracy and H. Widom, \emph{Level-spacing distributions and the
Airy kernel}, Commun. Math. Phys. \textbf{159} (1994), 151--174.

\bibitem[TW2]{TW2}
C. Tracy and H. Widom, \emph{On Orthogonal and Sympletic Matrix
Ensembles}, Communications in Mathematical Physics \textbf{177}
(1996), 727--754.

\bibitem[TW3]{TW3}
C. Tracy and H. Widom, \emph{Distribution functions for largest
eigenvalues and their applications}, ICM Vol. I (2002), 587--596.

\bibitem[VBAB]{VBAB}
M. G. Vavilov, P. W. Brouwer, V. Ambegaokar and C. W. J. Beenakker,
\emph{Universal gap fluctuations in the superconductor proximity
effect}, Phys. Rev. Letts. \textbf{86} (2001), 874--877.

\bibitem[Wig]{Wig}
E. Wigner, \emph{Statistical Properties of real symmetric matrices}.
Pages 174--184 in \emph{Canadian Mathematical Congress Proceedings},
University of Toronto Press, Toronto, 1957.

\end{thebibliography}
\end{document}